\documentclass[12pt,amsbook]{article}

\usepackage[english,brazilian]{babel}
\usepackage[latin1]{inputenc}
\usepackage[T1]{fontenc}
\usepackage[dvips]{graphics}
\usepackage{graphicx}
\usepackage{amssymb}
\usepackage{amsfonts}
\usepackage{amsmath}
\usepackage{amsthm}
\usepackage{verbatim}
\usepackage{epsfig}
\usepackage{setspace}
\usepackage{yhmath}         

\newtheorem{teo}{Theorem}[section]

\newtheorem{lem}[teo]{Lemma}
\theoremstyle{definition}
\newtheorem{defi}[teo]{Definition}

\newcommand{\N}{\mathbb{N}}
\newcommand{\R}{\mathbb{R}}

\newcommand{\Z}{{\cal{Z}}}
\newcommand{\M}{{\cal{M}}}

\newcommand{\dt}{\,dt}

\newcommand{\nl}{\mbox{}\\}

\begin{document}

%
%

\mbox{} \vspace{-2.000cm} \\
\begin{center}
\mbox{\Large \bf %
Some properties for the Steklov averages} \\
\nl
\mbox{} \vspace{-0.350cm} \\
{\large \sc J.\;Q.\;Chagas,}$\mbox{}^{\!\:\!1}$
{\large \sc N.\;M.\;L.\;Diehl}$\mbox{}^{\;\!2}$
{\large \sc and P.\;L.\;Guidolin}$\mbox{}^{\;\!3}$ \\
\mbox{} \vspace{-0.125cm} \\
$\mbox{}^{1}${\small
Departamento de Matem\'atica e Estat\'istica} \\
\mbox{} \vspace{-0.685cm} \\
{\small
Universidade Estadual de Ponta Grossa} \\
\mbox{} \vspace{-0.685cm} \\
{\small
Ponta Grossa, PR 84030-900, Brazil} \\
\mbox{} \vspace{-0.350cm} \\
$\mbox{}^{2}${\small
Instituto Federal de Educa\c{c}\~ao, Ci\^encia e Tecnologia do Rio Grande do Sul} \\
\mbox{} \vspace{-0.685cm} \\
{\small
Canoas, RS 92412-240, Brazil} \\
\mbox{} \vspace{-0.350cm} \\
$\mbox{}^{3}${\small
Instituto Federal de Educa\c{c}\~{a}o, Ci\^encia e Tecnologia do Rio Grande do Sul} \\
\mbox{} \vspace{-0.685cm} \\
{\small Viam\~{a}o, RS, 94410-970, Brazil} \\

%
%

\nl
\mbox{} \vspace{-0.400cm} \\
{\bf Abstract} \\
\mbox{} \vspace{-0.525cm} \\
\begin{minipage}[t]{11.750cm}
{\small
\mbox{} \hspace{+0.150cm}
We derive and present a collection of properties about the Steklov averages, including some results about the derivation with respect to spatial variables, and with respect to time, and a form of the fundamental theorem of the calculus. \\
}
\end{minipage}
\end{center}

%
%

\section{Introduction}

In this work we'll derive and present a collection of properties for the Steklov averages, which are an important regularization technique used currently in study of PDE's theory, but let us start by some brief notes about this mathematical tool and its proponent.

The Steklov average (or Steklov mean function) was introduced by V. A. Steklov in 1907 (see \cite{Stekloff1907}) for the study of the problem of expanding a given function into a series of eigenvalues defined by a 2nd-order ordinary differential operator; and its definition appears in $\S 67$ of \cite{Akhiezer1992} (along with some properties in $\S 83$). We reproduce this definition here:
\begin{center}
\hspace{4cm}
\begin{minipage}[t]{11.75cm}
{\small
Suppose the function $f(t)$, defined along the entire real axis, belongs to $L(a,b)$, for all finite values of $a$, $b$. Given any positive $h$, let us now construct the function
\begin{equation*}
 \displaystyle
 f_{h}(t) = \frac{1}{h} \int_{t-\frac{h}{2}}^{t+\frac{h}{2}} f(u) \, du
          = \frac{1}{h} \int_{-\frac{h}{2}}^{\frac{h}{2}} f(t+v) \, dv.
\end{equation*}
}
\end{minipage}
\end{center}

Vladimir Andreevich Steklov (1864 - 1926) was an out-standing russian mathematician who made many important contributions to mathematical physics (the Steklov average is only one of the mathematical notions associated with his name). Moreover, in 1921 Steklov founded the Physical-Mathematical Institute in Petrograd. Today, a famous institute of mathematics in Moscow has Steklov's name. More about life and work of V. A. Steklov can be seen in \cite{Kuznetsovetall2014}, \cite{Kuznetsov2014} and \cite{Sinai2003}.

In the present-day, Steklov averages are a very useful starting point for the derivation of a number of important solution properties for the PDE's (see \cite{DiBenedetto1993}, \cite{Urbano2008} or \cite{WuZhaoYinLi2001}).
In \cite{DiBenedetto1993}, the Steklov average is used to define a local weak solution that involves the time derivative $u_t$ for quasilinear degenerate or singular parabolic equations, widely used to derive solution properties;
in \cite{WuZhaoYinLi2001}, it's used in Chapter 2, which treats the non-newtonian filtration equations;
and in \cite{Urbano2008} it's used to define a weak solution for its central example, the parabolic p-Laplacian equation, because is proper to expose the idea of intrinsic scaling.
Recent uses of the Steklov Average in treatment of PDE's problems can be seen in \cite{ChagasGuidolinJZingano2017} and \cite{ChagasGuidolinZingano2017}.

However, only some properties of the Steklov averages are readily found in the literature. For this reason, in this work we've proposed to obtain and present a collection of important and useful operational properties for the Steklov averages.

Here is a brief description of what follows. In section 2, we present some convergence results for the Steklov average. In the section 3, we present the pointwise value form for the Steklov average. In section 4, we present some properties about the differentiability of the Steklov average, for spatial variables, and with respect to time $t$. Finally, in the section 5, we present some properties about the integration of the Steklov average, including forms of the fundamental theorem of calculus and integration by parts.

%
%

\section{About the convergence of Steklov averages}\label{2}

We'll start this section presenting the definition of the Steklov average of a function. In the sequence, we'll present results about the convergence of such avarages.

The space of the measurable sets in ${\R}^{n}$ it will be denoted by $\mathcal{M}({\R}^{n})$.

\begin{defi}\label{1} Let $[a,b]$ a compact interval in $\R$, $E \in \mathcal{M} (\R^n)$ and $1 \leq q \leq \infty$. Given $v(\cdot , t) \in C^0([a,b], L^q(E))$, we define (for each $0 < h < b-a$) the \textbf{Steklov average} $v_{h}(\cdot,t)$ of an function $v$ by
\begin{equation}\tag{2.1}
\displaystyle
v_{\mbox{}_h}(\cdot, t  ) = \frac{1}{h} \int_t^{t+h} v(\cdot ,s) \, ds, \quad \text{ for } \ \ a \leq t \leq b-h.
\end{equation}
\end{defi}

%
%

\vspace{0.5cm}
The first result is presented in:

\begin{lem}\label{2.2}
Given $ v(\cdot ,t) \in C^0([a,b], L^q(E))$, we have
\begin{equation*}
v_{\mbox{}_h}(\cdot,t) \in C^0([a,b-h],L^q(E))
\end{equation*}
(moreover, $v_{\mbox{}_h}(\cdot,t)$ is Lipschitz continous in the interval $[a,b-h]$),

\vspace{0.3cm} \noindent
and, for each $a \leq t < b$, we have
\begin{equation}\tag{2.2}\label{2.2}
\|v_{\mbox{}_h}(\cdot,t) - v(\cdot,t)\|_{L^q(E)} \rightarrow 0, \quad
\text{ when } h \rightarrow 0,
\end{equation}
uniformly in $t \in [a,b-\epsilon]$, for each $0 < \epsilon < h$.
\end{lem}

\vspace{0.3cm} \noindent
{\bf Proof.}
Let's prove this lemma in three assertions.

\vspace{0.3cm} \noindent
{\it Assertion $(i):$} $\, \displaystyle v_{\mbox{}_h}(\cdot, t  ) \in L^q(E), \ \forall \ t \in [a,b-h].$

\vspace{0.3cm} \noindent
Indeed, given $a \leq t \leq b-h $, we have
\begin{equation*}
\begin{aligned}
\displaystyle \| v_{\mbox{}_h}(\cdot,t)\|_{L^q(E)} & \leq \frac{1}{h} \left\| \int_t^{t+h}  v(\cdot,s) \, ds \right\|_{L^q(E)} \leq \frac{1}{h}  \int_t^{t+h} \|  v(\cdot,s) \|_{L^q(E)} \, ds\\
 \displaystyle & \leq \frac{1}{h}  \int_a^{b} \|  v(\cdot,s) \|_{L^q(E)}  \, ds \, \leq \, \infty,\\
\end{aligned}
\end{equation*}
so that
$\displaystyle \| v_{\mbox{}_h}(\cdot,t)\|_{L^q(E)} \leq M_h, \ \forall \ t \in [a,b-h] $, where $\displaystyle  M_h = \frac{1}{h} \int_a^{b} \|  v(\cdot,s) \|_{L^q(E)}  \, ds.$

\vspace{0.3cm} \noindent
{\it Assertion $(ii):$} $\, \displaystyle v_{\mbox{}_h}(\cdot,t) \in C^0([a,b-h],L^q(E)).$

\vspace{0.3cm} \noindent
Let $\displaystyle M := \max_{a \leq t \leq b} \|  v(\cdot,t) \|_{L^q(E)}$. Then, given $t_1 < t_2 \in [a,b-h]$, with $t_2 - t_1 \leq h$, we have
\begin{equation*}
\begin{aligned}
\displaystyle v_{\mbox{}_h}(\cdot,t_2) - \displaystyle v_{\mbox{}_h}(\cdot,t_1) & = \frac{1}{h} \int_{t_2}^{t_2+h} v(\cdot,t) \dt -\frac{1}{h} \int_{t_1}^{t_1+h} v(\cdot,t) \dt \\
& = \frac{1}{h} \int_{t_1+h}^{t_2+h} v(\cdot,t) \dt -\frac{1}{h} \int_{t_1}^{t_2}  v(\cdot,t) \dt,
\end{aligned}
\end{equation*}
and thus
\begin{equation*}
\begin{aligned}
\displaystyle \| v_{\mbox{}_h}(\cdot,t_2) - \displaystyle v_{\mbox{}_h}(\cdot,t_1) \|_{\mbox{}_{L^q(E)}} & \leq \frac{1}{h} \int_{t_2}^{t_2+h} \!\! \| v(\cdot,t) \|_{\mbox{}_{L^q(E)}} \dt  + \frac{1}{h} \int_{t_1}^{t_1+h} \!\! \| v(\cdot,t) \|_{\mbox{}_{L^q(E)}} \dt \\
&\leq \frac{1}{h} \int_{t_1+h}^{t_2+h} \!\! M \dt + \frac{1}{h} \int_{t_1}^{t_2} \!\! M \dt = \frac{2}{h}M \,|t_2 - t_1|.
\end{aligned}
\end{equation*}

\vspace{0.3cm} \noindent
{\it Assertion $(iii):$} $\displaystyle v_{\mbox{}_h}(\cdot,t) \rightarrow v(\cdot,t)$ \ in $L^q(E)$ as $h \rightarrow 0$, for each $t \in [a,b)$.

\vspace{0.3cm} \noindent
Indeed, for $T \in [a,b)$, we define $\varepsilon_T := b - T > 0$ and take $0< h < \varepsilon_T$. Given $\varepsilon > 0$, let $\delta >0$ such that
$$
\| v(\cdot,t) - v(\cdot,s) \|_{L^q(E)} \leq \varepsilon, \quad \forall \ s,t \in [a,b] \ \text{ with } |s-t| \leq \delta \leq \varepsilon_T.
$$
Then, for each $0 < h \leq \delta$ and $\forall \ t \in [a,T]$, we obtain
\begin{equation*}
\begin{aligned}
\displaystyle \| v_{\mbox{}_h}(\cdot,t) - v(\cdot,t) \|_{L^q(E)} & \leq \, \left\| \frac{1}{h} \int_{t}^{t+h} v(\cdot,s) \, ds - \frac{1}{h} \int_{t}^{t+h} v(\cdot,t) \, ds  \right\|_{L^q(E)} \\
& \leq \, \frac{1}{h} \int_{t}^{t+h} \| v(\cdot,s) - v(\cdot,t)\|_{L^q(E)} \, ds \\
& \leq \, \frac{1}{h} \int_{t}^{t+h} \varepsilon \, ds = \varepsilon ,
\end{aligned}
\end{equation*}
i.e., for all $t \in [a,T]$ we have that
$$
\| v_{\mbox{}_h}(\cdot,t) - v(\cdot,t) \|_{L^q(E)} \leq \varepsilon, \quad \forall \ 0 < h \leq \delta.
$$
\vspace{-1.25cm}
\begin{flushright}
(Lema 2.2) \ $\square$
\end{flushright}

For the proof of the next lemma, it's convenient define $v_{h}(\cdot,t)$ as follows:

\vspace{0.3cm} \noindent
\begin{defi}\label{3}
Let $I \subset \R$ any interval, $E \in \mathcal{M}({\R}^{n})$, $1 \leq q,r \leq \infty$ and $h> 0$. Given $v(\cdot,t) \in L^r(I, L^q(E))$ we define $v_h(\cdot,t) \in C^0(I,L^q(E) )$ by
\begin{equation}\tag{2.3$a$}\label{2.3a}
 \displaystyle
 v_h(\cdot,t) = \frac{1}{h} \int^{t+h}_h \tilde{v}(\cdot,s) \, ds, \quad \text{for each} \ a \leq t \leq b-h,
\end{equation}
where $\tilde{v}(\cdot,t) \in L^r(\R, L^q(E))$ is defined by
\begin{equation}\tag{2.3$b$}\label{2.3b}
\tilde{v}(\cdot,t)=
\left\{
 \begin{aligned}
  v(\cdot,t), & \ \text{ se } t \in I,  \\
  0, & \ \text{ se } t \in \R \! \setminus \! I.
 \end{aligned}
\right.
\end{equation}
\end{defi}

%
%

\vspace{0.3cm}
\begin{lem}\label{2.4}
Given an interval $I \subseteq \R$, $E \in \mathcal{M}(\R^n)$, $1 \leq q ,r \leq \infty$ and $h> 0$, let $v(\cdot,t) \in L^r(I,L^q(E))$. Then, $v_h(\cdot,t)$ as defined in (\ref{2.3a})-(\ref{2.3b}) satisfies:
\begin{equation}\tag{2.4$a$}\label{2.4a}
 v_h(\cdot,t) \in L^q(E), \ \forall \ t \in I,
\end{equation}
with
\begin{equation*}
 \|v_h(\cdot,t)\|_{L^q(E)} \leq \frac{1}{h^{1/r}} \|v(\cdot,t)\|_{L^r(I,L^q(E))}, \text{ if } r \neq \infty,
\end{equation*}
\begin{equation*}
 \|v_h(\cdot,t)\|_{L^q(E)} \leq  \|v(\cdot,t)\|_{L^{\infty}(I,L^q(E))}, \text{ if } r = \infty,
\end{equation*}
\begin{equation}\tag{2.4$b$}\label{2.4b}
 v_h(\cdot,t) \in C^0(I,L^q(E)) \cap L^{\infty}(I,L^q(E)),
\end{equation}
\begin{equation}\tag{2.4$b$'}\label{2.4b'}
 v_h(\cdot,t) : I \to L^q(E) \ \text{ is uniformly continuous on } I.
\end{equation}
\begin{center}
(and, if $r = \infty$, then $v_h(\cdot,t) : I \to L^{q}(E)$ is Lipschitz continuous on $I$),
\end{center}
\begin{equation}\tag{2.4$c$}\label{2.4c}
 v_h(\cdot,t) \in L^r(I,L^q(E)),
\end{equation}
and
\begin{equation*}
 \|v_h\|_{L^r(I,L^q(E))} \leq \|v\|_{L^r(I,L^q(E))}.
\end{equation*}
\end{lem}

\vspace{0.3cm} \noindent
{\bf Proof.}

\vspace{0.3cm} \noindent
{\it Proof of (\ref{2.4a}): } If $r = \infty$, for each $t \in I$ we have
\begin{equation*}
\begin{aligned}
\displaystyle
{\| v_h(\cdot,t) \|}_{L^q(E)}
& \leq \frac{1}{h} \int_{t}^{t+h} {\|\tilde{v}(\cdot,s)\|}_{L^q(E)} \, ds \\
& \leq \frac{1}{h} \int_{t}^{t+h} {\|{v}\|}_{L^{\infty}(I,L^q(E))} \, ds
= {\|{v}\|}_{L^{\infty}(I,L^q(E))}.
\end{aligned}
\end{equation*}
If $1 \leq r < \infty$, for each $t \in I$ we have
\begin{equation*}
\begin{aligned}
 \displaystyle \| v_h (\cdot,t)\|_{L^q(E)} & \leq \frac{1}{h} \int_t^{t+h} \|\tilde{v}(\cdot,s)\|_{L^q(E)} \, ds \leq \frac{1}{h^{1-1/r'}}\left( \int_t^{t+h} {\|\tilde{v}(\cdot,s)\|}^{r}_{L^q(E)} \, ds \right)^{1/r} \qquad  \\
 & = \frac{1}{h^{1/r}}\!\! \left( \int_t^{t+h} {\|\tilde{v}(\cdot,s)\|}^{r}_{L^q(E)} \, ds \right)^{\!\! 1/r} \!\!\!\! \leq \!  \frac{1}{h^{1/r}}\!\!\left( \int_t^{t+h} {\|{v}(\cdot,s)\|}^{r}_{L^q(E)} \, ds \right)^{\!\! 1/r}
\end{aligned}
\end{equation*}
i.e.,
$$
\displaystyle \| v_h (\cdot ,t )\|_{L^q(E)} \leq  \frac{1}{h^{1/r}} \|{v}\|_{L^r(I,L^q(E))} \quad \forall t \in I.
$$
In particular, for $1 \leq r  \leq \infty $ we have $ v_h (\cdot ,t ) \in L^{\infty}(I,L^q(E))$.

\vspace{0.3cm} \noindent
 {\it Proof of (\ref{2.4b}) and (\ref{2.4b'}): } If $r= \infty$, for each $t_1 < t_2 \in I$ with $|t_2 - t_1| < h$, we have
\begin{equation*}
\begin{aligned}
 \displaystyle
 v_h(\cdot,t_2) - v_h(\cdot,t_1) & = \frac{1}{h} \int_{t_2}^{t_2+h} \tilde{v}(\cdot,s) \, ds - \frac{1}{h} \int_{t_1}^{t_1+h} \tilde{v}(\cdot,s) \, ds \\
 & = \frac{1}{h} \int_{t_1+h}^{t_2+h} \tilde{v}(\cdot,s) \, ds - \frac{1}{h} \int_{t_1}^{t_2} \tilde{v}(\cdot,s) \, ds
\end{aligned}
\end{equation*}
therefore
\begin{equation*}
\begin{aligned}
 \displaystyle
 {\| v_h(\cdot,t_2) - v_h(\cdot,t_1)\|}_{L^q(E)} & \leq \frac{1}{h} \int_{t_1+h}^{t_2+h} {\|\tilde{v}(\cdot,s)\|}_{L^q(E)} \, ds + \frac{1}{h} \int_{t_1}^{t_2} {\|\tilde{v}(\cdot,s)\|}_{L^q(E)} \, ds \\
 & \leq \frac{1}{h} \int_{t_1+h}^{t_2+h} M \, ds + \frac{1}{h} \int_{t_1}^{t_2} M \, ds = \frac{2}{h} M |t_2 - t_1|,
\end{aligned}
\end{equation*}
where $M = \| v \|_{L^{\infty}(I,L^q(E))}$.
This shows that
\begin{equation*}
\displaystyle \| v_h(\cdot,t) - v_h(\cdot,s)\|_{L^q(E)} \leq \frac{2}{h}M|t-s|, \ \forall \ s,t \in I.
\end{equation*}

If $1\leq r < \infty$, for any $t_1 < t_2 \in I$ with $|t_2 - t_1| \leq h$, we have that
\begin{equation*}
\begin{aligned}
 \displaystyle
 v_h(\cdot,t_2) - v_h(\cdot,t_1) & = \frac{1}{h} \int_{t_2}^{t_2+h} \tilde{v}(\cdot,s) \, ds - \frac{1}{h} \int_{t_1}^{t_1+h} \tilde{v}(\cdot,s) \, ds \\
 & = \frac{1}{h} \int_{t_1+h}^{t_2+h} \tilde{v}(\cdot,s) \, ds - \frac{1}{h} \int_{t_1}^{t_2} \tilde{v}(\cdot,s) \, ds
 \end{aligned}
\end{equation*}
as before, but now follows that
\begin{equation*}
\begin{aligned}
 \displaystyle
 {\| v_h(\cdot,t_2) - v_h(\cdot,t_1)\|}_{L^q(E)} \leq \frac{1}{h} \int_{t_1+h}^{t_2+h} {\|\tilde{v}(\cdot,s)\|}_{L^q(E)} \, ds + \frac{1}{h} \int_{t_1}^{t_2} {\|\tilde{v}(\cdot,s)\|}_{L^q(E)} \, ds \\
 \leq \frac{1}{h^{1/r}} {\left( \int_{t_1+h}^{t_2+h} {\|\tilde{v}(\cdot,s)\|}^{r}_{L^q(E)} \, ds \right)}^{\!\!1/r}\!\!\!\! + \frac{1}{h^{1/r}} {\left(\int_{t_1}^{t_2} {\|\tilde{v}(\cdot,s)\|}^{r}_{L^q(E)} \, ds \right)}^{\!\!1/r}.
\end{aligned}
\end{equation*}

As $\displaystyle {\|\tilde{v}(\cdot,s)\|}_{L^q(E)}^{r} \in L^1(\R)$, given $\varepsilon > 0 $, there exists $\delta > 0 $ such that
\begin{equation*}
 \displaystyle
 \int_J {\| \tilde{v}(\cdot,s)\|}_{L^q(E)}^{r} \leq {\varepsilon}^{r}
\end{equation*}
\vspace{-0.2cm}
whenever $J \in \mathcal{M}(\R) $ and $|J| \leq \delta$.

Therefore, for any $t_1,t_2 \in I$ with $|t_2 - t_1| \leq \delta$, we obtain
\begin{equation*}
\begin{aligned}
 \displaystyle
 {\| v_h(\cdot,t_2) - v_h(\cdot,t_1)\|}_{L^q(E)} \leq \frac{1}{h^{1/r}} {({\varepsilon}^{r})}^{1/r} + \frac{1}{h^{1/r}} {({\varepsilon}^{r})}^{1/r} = \frac{2}{h^{1/r}} \varepsilon.
\end{aligned}
\end{equation*}
This shows that $v_h(\cdot, t) : I \rightarrow L^q(E)$ is uniformly continuous in $I$.

\vspace{0.2cm}
As $v_h(\cdot, t) : I \rightarrow L^q(E)$ is also limited in $I$, for each $ 1 \leq q,r \leq  \infty $, we have that
$v_h(\cdot, t) : I \rightarrow L^q(E)$ is limited and uniformly continuous in $I$; and for $1 \leq q \leq \infty , \, r = \infty$, we have that
$v_h(\cdot, t) : I \rightarrow L^q(E)$ is limited and globally Lipschitz in $I$.

\vspace{0.3cm} \noindent
{\it Proof of (\ref{2.4c}):} The case where $r= \infty$ already been shown in (\ref{2.4a}). Consider then $1 \leq r < \infty$.
From (\ref{2.3a}) we obtain, $\forall \ t \in I$, that
\begin{equation*}
 \displaystyle
 {\|v_h(\cdot,t)\|}_{L^q(E)} \leq \frac{1}{h} \int_{t}^{t+h} {\| \widetilde{v}(\cdot,s) \|}_{L^q(E)} \, ds \leq \frac{1}{h^{1/r}} {\left( \int^{t+h}_t {\| \tilde{v}(\cdot,s) \|}_{L^q(E)}^{r} \, ds \right)}^{1/r}.
\end{equation*}
Therefore
\begin{equation*}
 \displaystyle
 {\|v_h(\cdot,t)\|}_{L^q(E)}^{r} \leq \frac{1}{h} \left( \int^{t+h}_t {\| \tilde{v}(\cdot,s) \|}_{L^q(E)}^{r} \, ds \right), \quad \forall \ t \in I,
\end{equation*}
and thus, if $I = [a,b], \,(a,b], \, [a,b)$, or $(a,b)$, for $- \infty < a < b < \infty$, follows that

\begin{equation*}
\begin{aligned}
 \displaystyle
 \int_I {\| v_h(\cdot,t) \|}_{L^q(E)}^{r} \dt
 & = \frac{1}{h} \int_a^b \left( \int^{t+h}_t {\| \tilde{v}(\cdot,t) \|}_{L^q(E)}^{r} \, ds \right) \dt \\
 & = \frac{1}{h} \int_a^b \bigl( V(t+h) - V(t) \bigr) \dt \\
\end{aligned}
\end{equation*}
\begin{equation*}
\begin{aligned}
 \qquad \qquad \qquad \qquad
 & = \frac{1}{h} \int_{a+h}^{b+h} V(t) \dt - \frac{1}{h} \int_a^b V(t) \dt \\
 & \leq \frac{1}{h} \int_{b}^{b+h} V(t) \dt = V(b) = \int_a^b \|v(\cdot,s)\|^r_{L^q(E)} \, ds,
\end{aligned}
\end{equation*}
where $\displaystyle V(t) := \int_a^t \| \tilde{v}(\cdot,s)\|^r_{L^q(E)} \, ds , \ \forall \ t \geq a$; i.e., when $I$ is bounded, we have
\begin{equation}\tag{2.4$d$}\label{2.4d}
 \displaystyle
  \int_I {\|v_h(\cdot,s)\|}_{L^q(E)}^{r} \dt \leq \int_I {\|v(\cdot,s)\|}_{L^q(E)}^{r} \, ds.
\end{equation}

We'll now extend (\ref{2.4d}) for the cases $(-\infty,b)$, $(- \infty ,b]$; and $(a,\infty)$, $[a, \infty)$; with $a,b \in \R$.

In the cases $(-\infty,b)$ and $(- \infty ,b]$, we obtain
\begin{equation*}
\begin{aligned}
 \displaystyle \
 \int_I {\|v_h(\cdot,t)\|}_{L^q(E)}^{r} \dt & = \lim_{a \rightarrow - \infty} \int_a^b {\|v_h(\cdot,t)\|}_{L^q(E)}^{r} \dt \leq \lim_{a \rightarrow - \infty} \int_a^b {\|v(\cdot,t)\|}_{L^q(E)}^{r} \dt \\
 & = \int_{-\infty}^b {\|v(\cdot,t)\|}_{L^q(E)}^{r} \dt = \int_I {\|v(\cdot,t)\|}_{L^q(E)}^{r} \dt,
\end{aligned}
\end{equation*}
and in the cases $(a,\infty)$ and $[a,\infty)$ we obtain
\begin{equation*}
\begin{aligned}
 \displaystyle
 \int_I {\|v_h(\cdot,t)\|}_{L^q(E)}^{r} \dt & = \lim_{b \rightarrow \infty} \int_a^b {\|v_h(\cdot,t)\|}_{L^q(E)}^{r} \dt \leq \lim_{b \rightarrow \infty} \int_a^{b+h} {\|v(\cdot,t)\|}_{L^q(E)}^{r} \dt \\
 & = \int_a^{\infty} {\|v(\cdot,t)\|}_{L^q(E)}^{r} \dt = \int_I {\|v(\cdot,t)\|}_{L^q(E)}^{r} \dt,
\end{aligned}
\end{equation*}
where the inequality is obtained using $V(t)$, again, and the fact that
\begin{equation*}
 \displaystyle
 \frac{1}{h} \int_b^{b+h} V(t+h) \dt \leq \frac{1}{h} \int_b^{b+h} V(b+h) \dt = V(b+h) = \int_a^{b+h} {\|v(\cdot,t)\|}_{L^q(E)}^{r} \dt.
\end{equation*}

Finally, for the case $I= \R$, we have
\begin{equation*}
\begin{aligned}
 \displaystyle
 \int_I {\|v_h(\cdot,t)\|}_{L^q(E)}^{r} \dt & = \lim_{\stackrel{a \rightarrow - \infty }{b \rightarrow  \infty}} \int_a^b {\|v_h(\cdot,t)\|}_{L^q(E)}^{r} \dt \leq \lim_{\stackrel{a \rightarrow - \infty }{b \rightarrow \infty}} \int_a^{b+h} {\|v(\cdot,t)\|}_{L^q(E)}^{r} \dt \\
 & = \int_{-\infty}^{\infty} {\|v(\cdot,t)\|}_{L^q(E)}^{r} \dt = \int_I {\|v(\cdot,t)\|}_{L^q(E)}^{r} \dt.
\end{aligned}
\end{equation*}
Thus, for each interval $I \subset \R$, we have
\begin{equation*}
 \displaystyle
 \int_I {\|v_h(\cdot,t)\|}_{L^q(E)}^{r} \dt \leq \int_I {\|v(\cdot,t)\|}_{L^q(E)}^{r} \dt.
\end{equation*}

\vspace{-0.75cm}
\begin{flushright}
(Lemma \ref{2.4}) \ $\square$
\end{flushright}

%
%

The proof of the next lemma requires $1 \leq r < \infty$.

\vspace{0.1cm}
\begin{lem}\label{2.5}
Given any interval $I \subset \R$, $E \in \mathcal{M}(\R^n), \, 1 \leq q \leq \infty, \,$ $ 1 \leq r < \infty$ and $ h> 0 $, take $v(\cdot,t) \in L^r(I, L^q(E))$ and consider $v_h(\cdot,t) \in C^0(I,L^q(E)) \cap L^{\infty}(I,L^q(E))$ (see Lemma \ref{2.4}).

Then, when $h \to 0$ we have
\vspace{-0.5cm}
\begin{equation}\tag{2.5$a$}\label{2.5a}
 v_h \rightarrow v \ \text{ in } L^r(I,L^q(E)),
\end{equation}

\vspace{-0.5cm} \noindent
and

\vspace{-0.3cm}
\begin{equation}\tag{2.5$a$'}\label{2.5a'}
 \displaystyle
 \int_I {\|v_h(\cdot,t) - v(\cdot, t)\|}_{L^q(E)}^{r} \dt \rightarrow 0.
\end{equation}
\end{lem}

\vspace{0.3cm} \noindent
{\bf Proof.} Let $\ring{I}$ the interior of the set $I$.

Initially, let's assume that $v(\cdot,t) \in C^0_c(\ring{I}, L^q(E))$, that is, $v(\cdot,t) \in C^0(I,L^q(E))$, with $v(\cdot,t) =0 \ \ \forall  \ t \in I \! \setminus \! [a,b]$, for some compact interval $[a,b] \subset \ring{I}$.

Then, we have that
\begin{equation*}
 \displaystyle
 {\|v_h(\cdot,t) - v(\cdot,t) \|}_{L^q(E)} \rightarrow 0
\end{equation*}
uniformly in $t \in I$ as $h \rightarrow 0$, as has been proved in the item (\ref{2.2}) of Lemma \ref{2.2}.

\vspace{0.3cm}
Indeed, let a compact interval $[a,b] \subset \ring{I}$ containing the support of $v(\cdot,t)$, and take a compact interval $[\alpha , \beta] \subset \ring{I}$ with $\alpha < a$ and $\beta > b$. For each $h>0$ with
\begin{equation*}
 \displaystyle
 h \leq \min \{ a-\alpha,\beta - b \} := \bar{h},
\end{equation*}
we have that $v_h(\cdot,t) = 0 = v(\cdot,t)$ if $t \in I$ satisfies $t \leq \alpha$ or $t \geq \beta$, and therefore,
\begin{equation*}
 \displaystyle
 {\|v_h(\cdot,t) - v(\cdot,t)\|}_{L^q(E)} = 0 \ \ \forall \ 0 < h < \bar{h},
\end{equation*}
for each $t \in I$ with $t \leq \alpha$ or $t \geq \beta$.

On the other hand, for $t \in [\alpha , \beta]$, we may proceed as follows: given $\varepsilon > 0$, let $\delta >0$ be small enough for that
\begin{equation*}
 \displaystyle
 {\|v(\cdot,s) - v(\cdot,r) \|} \leq \varepsilon \ \  \forall \ s,r \in [\alpha , \beta]
\end{equation*}
with $|s-r| \leq \delta$ and $\delta \leq \bar{h}$. Then, for every $t \in [\alpha , \beta]$ and $\forall \ 0 < h < \delta \leq \bar{h}$, we have that
\begin{equation*}
 \displaystyle
 {\|v(\cdot,s) - v(\cdot,r) \|}_{L^q(E)} \leq \frac{1}{h} \int_t^{t+h} \! {\|v(\cdot,s) - v(\cdot,r) \|}_{L^q(E)} \, ds \leq \frac{1}{h} \int_t^{t+h} \! \varepsilon \, ds = \varepsilon.
\end{equation*}
For $v(\cdot,t) \in C_c^0(\ring{I},L^q(E))$ this shows that
\begin{equation}\tag{2.5$b$}\label{2.5b}
 \displaystyle
 \int_I {\|v_h(\cdot,t) - v(\cdot,t)\|}_{\mbox{}_{L^q(E)}}^{r} = \int_{\alpha}^{\beta} {\|v_h(\cdot,t) - v(\cdot,t)\|}_{\mbox{}_{L^q(E)}}^{r} \rightarrow 0 \ \text{ when } h \rightarrow 0,
\end{equation}
which results in (\ref{2.5a}) e (\ref{2.5a'}), in the case where $v(\cdot,t) \in C^{0}_{c}(\ring{I},L^{q}(E))$.

\vspace{0.3cm}
In the general case $v(\cdot,t) \in L^r(I, L^q(E))$, we may proceed as follows: given $\varepsilon > 0 $, we can take $w(\cdot,t) \in C_c^0(\ring{I},L^q(E))$ such that
\begin{equation*}
 \displaystyle
 \int_I {\| v(\cdot,t) - w(\cdot,t) \|}^{r}_{L^q(E)} \leq {\left( \, \frac{\varepsilon}{3} \, \right)}^r,
\end{equation*}
(this follows because $C_c^{\infty}(\ring{I}, L^q(E)) $ is dense in $L^r(I, L^q(E)), \  \forall \ 1 \leq r < \infty$), i.e.,
\begin{equation*}
 \displaystyle
 { \| v - w \| }_{L^{r}(I,L^{q}(E))} \leq \frac{\varepsilon}{3}.
\end{equation*}
In particular, by (\ref{2.4c}) of Lemma \ref{2.4}, we also have that ${ \| v_h - w_h \|}_{L^{r}(I,L^{q}(E))} \leq \frac{\varepsilon}{3}$.

Therefore, we have
\begin{equation*}
\begin{aligned}
 \displaystyle
 {\| v_h - v \|}_{L^{r}(I,L^{q}(E))} & \leq {\| v_h - w_h \|}_{L^{r}(I,L^{q}(E))} + {\| w_h - w \|}_{L^{r}(I,L^{q}(E))} + {\| w - v \|}_{L^{r}(I,L^{q}(E))} \\
  &  \leq \frac{\varepsilon}{3} +  {\| w_h - w \|}_{L^{r}(I,L^{q}(E))} + \frac{\varepsilon}{3}.
\end{aligned}
\end{equation*}
As $w(\cdot,t) \in C^0_c(\ring{I},L^q(E))$, from (\ref{2.5b}), we know that ${\| w_h - w \|}_{L^{r}(I,L^{q}(E))} \rightarrow 0 $ as $h \rightarrow 0$. Therefore, taking $h_0 > 0$ small enough for occurs
\begin{equation*}
 \displaystyle
 {\| w_h - w \|}_{L^{r}(I,L^{q}(E))} \leq \frac{\varepsilon}{3} \quad \forall  \ 0 < h \leq h_0,
\end{equation*}
we obtain that
$\displaystyle {\| v_h - v\|}_{L^{r}(I,L^{q}(E))} \leq \varepsilon, \ \ \forall \ 0 < h \leq h_0.$

\vspace{0.2cm}
\vspace{-0.75cm}
\begin{flushright}
(Lemma $2.5$) \ $\square$
\end{flushright}

%
%

Finally, we can show the next Lemma (since $1 \leq r < \infty$), which has the harder proof to obtain among the presented results until here.

\begin{lem}\label{2.6}
Given $I \subset \R$ (any interval), $E \in \mathcal{M}(\R^n)$, $1 \leq q \leq \infty$, $1 \leq r < \infty$, $h>0$, and $v(\cdot,t) \in L^r(I,L^q(E))$ arbitrary, consider $v_h \in C^0(I,L^q(E)) \cap L^{\infty}(I,L^q(E)) $ defined as in (\ref{2.3a}).

Then, there exists $Z \subseteq I$ with zero measure such that, for each $t \in I \! \setminus \! Z$, we have
\begin{equation*}
 \displaystyle
 v(\cdot,t) \in L^q(E), \ \text{  and }
\end{equation*}
\begin{equation*}
 \displaystyle
 {\| v_h(\cdot,t) - v(\cdot,t)\|}_{\mbox{}_{L^q(E)}} \rightarrow 0 \ \text{  as } h \rightarrow 0.
\end{equation*}
\end{lem}

\vspace{0.3cm} \noindent
{\bf Proof.}
Given $v(\cdot,t) \in L^r(I,L^q(E))$, with $1 \leq r < \infty$, we take a sequence of smooth approximations $w_m(\cdot,t) \in C_c^0(\ring{I},L^q(E)), \ \forall \ m \in \N$, such that
\begin{equation}\tag{2.6$a$}\label{2.6a}
 \displaystyle
 {\| w_m - v \|}_{\mbox{}_{L^r(I,L^q(E))}} \rightarrow 0, \ \text{  as } m \rightarrow  \infty,
\end{equation}
and (passing to a subsequence, if necessary)
\begin{equation}\tag{2.6$a$'}\label{2.6a'}
 \displaystyle
 {\| w_m(\cdot,t)  - v(\cdot,t) \|}_{\mbox{}_{L^q(E)}} \rightarrow 0, \ \text{  as } m \rightarrow  \infty,
\end{equation}
for each $t \in I \! \setminus \! Z_0$, with $Z_0 \subseteq I$ of zero measure.

Observe that, because $v(\cdot,t) \in L^r(I,L^q(E))$, there exists $Z_v \subseteq I$ with zero measure such that $v(\cdot,t ) \in L^q(E), \ \forall \ t \in I \! \setminus \! Z_v$. The null set $Z_0$ in (\ref{2.6a'}) satisfies, in particular, $Z_v \subseteq Z_0.$

For each $m \geq 1$, by (\ref{2.6a}) we have that $\displaystyle {\| w_m(\cdot,t) - v(\cdot,t) \|}_{\mbox{}_{L^q(E)}} \in L^r(I)$, so that, by H\"{o}lder's inequality, we have
\begin{equation*}
 \displaystyle
 {\| w_m(\cdot,t) - v(\cdot,t) \|}_{\mbox{}_{L^q(E)}} \in L^1_{loc}(I).
\end{equation*}

By Lebesgue's differentiation theorem, there exists $Z_m \subseteq I$, with $|Z_m|=0$, and $Z_v \subseteq Z_m$ such that
\begin{equation}\tag{2.6$b$}\label{2.6b}
 \displaystyle
 \lim_{m \rightarrow 0} \frac{1}{h} \int_{t}^{t+h} {\| w_m(\cdot,s) - v(\cdot,s) \|}_{\mbox{}_{L^q(E)}} \, ds = {\| w_m(\cdot,t) - v(\cdot,t) \|}_{\mbox{}_{L^q(E)}}, \ \forall \ t \in I \! \setminus \! Z_m.
\end{equation}

Let then $\displaystyle  Z := Z_0 \cup \left( \bigcup_{m=1}^{\infty} Z_m  \right)$. Thus, we have $Z_v \subseteq  Z \subseteq I$, with $Z$ having zero measure.

We claim that, for each $t \in I \! \setminus \! Z$, we have
\begin{equation}\tag{2.6$c$}\label{2.6c}
 \displaystyle
 \lim_{m \rightarrow 0} {\| v_h(\cdot,t) - v(\cdot,t) \|}_{L^q(E)} = 0.
\end{equation}

Indeed, given $\hat{t} \in I \! \setminus \! Z$ and $\varepsilon >0$, we may proceed as follows: take $m_0 \geq 1$ big enough so that
\begin{equation}\tag{2.6$d$}\label{2.6d}
 \displaystyle
  {\| w_{m_0}(\cdot,\widehat{t}\,) - v(\cdot,\widehat{t}\,) \|}_{\mbox{}_{L^q(E)}}
  \leq \frac{\varepsilon}{4}.
\end{equation}
Let $w(\cdot,t) \in C^0_c(\ring{I},L^q(E))$ be given by $w(\cdot,t) = w_{m_0}(\cdot,t) \ \forall \ t \in I$. By (\ref{2.6b}), since $\hat{t} \in I \! \setminus \! Z_{m_0}$, we have
\begin{equation*}
 \displaystyle
 \frac{1}{h} \int_{\hat{t}}^{\hat{t}+h} {\| w(\cdot,s) - v(\cdot,s) \|}_{\mbox{}_{L^q(E)}} \, ds \rightarrow {\| w(\cdot,\hat{t}\,) - v(\cdot,\hat{t}\,) \|}_{\mbox{}_{L^q(E)}} \leq \frac{\varepsilon}{4}.
\end{equation*}
Hence, there exists $h_{\varepsilon} >0$ (by (\ref{2.6d})) small enough that we have
\begin{equation*}
 \displaystyle
 \frac{1}{h} \int_{\hat{t}}^{\hat{t}+h} {\| w(\cdot,s) - v(\cdot,s) \|}_{\mbox{}_{L^q(E)}} \, ds \leq \frac{\varepsilon}{2}, \ \ \forall \ 0 < h \leq h_{\varepsilon}.
\end{equation*}
Therefore, for any $0 < h \leq h_{\varepsilon}$, this gives
\begin{equation*}
 \begin{aligned}
  {\| w_h(\cdot,\hat{t}) - v_h(\cdot,\hat{t}) \|}_{\mbox{}_{L^q(E)}}
  & = \frac{1}{h} {\left\| \int_{\hat{t}}^{\hat{t}+h} \bigl( w(\cdot,s) - v(\cdot,s) \bigr) \, ds \right\|}_{\mbox{}_{L^q(E)}} \\
  & \leq \frac{1}{h} \int_{\hat{t}}^{\hat{t}+h} {\| w(\cdot,s) - v(\cdot,s) \|}_{\mbox{}_{L^q(E)}} \, ds \leq \frac{\varepsilon}{2}.
 \end{aligned}
\end{equation*}
Therefore, for all $h \in \ ]0,h_{\varepsilon}]$, we have
\begin{equation}\tag{2.6$e$}\label{2.6e}
 \begin{aligned}
  {\| v_h(\cdot,\hat{t}\,) - v(\cdot,\hat{t}\,) \|}_{L^q(E)}
  & \leq {\| v_h(\cdot,\hat{t}\,) - w_h(\cdot,\hat{t}\,) \|}_{L^q(E)}
  + {\| w_h(\cdot,\hat{t}\,) - w(\cdot,\hat{t}\,) \|}_{L^q(E)} \\
  & \qquad \qquad \qquad \qquad \qquad \qquad \quad + {\| w(\cdot,\hat{t}\,) - v(\cdot,\hat{t}\,) \|}_{L^q(E)} \\
  & \leq \frac{\varepsilon}{2} + {\| w_h(\cdot,\hat{t}\,) - w(\cdot,\hat{t}\,) \|}_{L^q(E)} + \frac{\varepsilon}{4}.
 \end{aligned}
\end{equation}

Because $w(\cdot,{t}) \in C^0_c(\ring{I},L^q(E))$, we clearly have
\begin{equation}\tag{2.6$f$}\label{2.6f}
 \displaystyle
 \lim_{h \rightarrow 0} {\| w_h(\cdot,t) - w(\cdot,t) \|}_{\mbox{}_{L^q(E)}} = 0,
 \ \text{ for each } t \in I,
\end{equation}
and hence, there exists $h_{\varepsilon \varepsilon} << 1$ such that
\begin{equation}\tag{2.6$f$\,'}\label{2.6f'}
 \displaystyle
 {\| w_h(\cdot,\hat{t}\,) - w(\cdot,\hat{t}\,) \|}_{\mbox{}_{L^q(E)}} \leq \frac{\varepsilon}{4}, \ \ \forall \ 0 < h < h_{\varepsilon \varepsilon}.
\end{equation}
From (\ref{2.6e}), (\ref{2.6f'}),
$\forall \ 0 < h \leq \min \{ h_{\varepsilon},h_{\varepsilon \varepsilon} \}$, we get
\begin{equation*}
 \begin{aligned}
  {\| v_h(\cdot,\hat{t}\,) - v(\cdot,\hat{t}\,) \|}_{\mbox{}_{L^q(E)}}
  & \leq \frac{\varepsilon}{2} + {\| w_h(\cdot,\hat{t}\,) - w(\cdot,\hat{t}\,) \|}_{\mbox{}_{L^q(E)}} + \frac{\varepsilon}{4} \\
  & \leq \frac{\varepsilon}{2} + \frac{\varepsilon}{4} + \frac{\varepsilon}{4} = \varepsilon \ \ \
  \forall \ 0 < h \leq \min \{ h_{\varepsilon},h_{\varepsilon \varepsilon}\}.
 \end{aligned}
\end{equation*}

This shows (\ref{2.6c}) for $t = \hat{t}$ (with $\hat{t} \in I \! \setminus \! Z$ arbitrary), as claimed.

\vspace{-0.5cm}
\begin{flushright}
(Lemma $2.6$) \ $\square$
\end{flushright}

%
%

\section{Pointwise values of the Steklov averages}\label{3}

\vspace{0.3cm}
Regarding the {\bf pointwise} values of $v_{h}(\cdot,t) \in C^{0}(I,L^{q}(E)) \cap L^{\infty}(I,L^{q}(E))$ (for a given $v(\cdot,t) \in L^{r}(I,L^{q}(E))$, where $1 \leq q,r \leq \infty$), we have, by Fubini's theorem, that the following result holds:

%
%

\begin{lem}\label{3.1}
For each $J \subseteq I$, where $J$ is a bounded interval, there exists $Z \subseteq E$, with $|Z|=0$ ($Z$ depending on $J$), such that
\begin{equation}\tag{3.1$a$}\label{3.1a}
 \displaystyle
 \int_{J} |v(x,t)| \, dt < \infty, \quad \forall  \ x \in E \! \setminus \! Z.
\end{equation}
It follows that there exists ${Z}_{*} \subseteq E$, with $|{Z}_{*}| = 0$, such that
\begin{equation}\tag{3.1$b$}\label{3.1b}
 \displaystyle
 \int_{t}^{t+h} {|\tilde{v}(x,t)|} \, ds < \infty, \ \ \forall \ x \in E \! \setminus \! {Z}_{*}, \ \ \forall \ t \in I, \ \ \forall \ h > 0.
\end{equation}
In particular, we have that the pointwise values of $v_{h}(\cdot,t)$ are given by
\begin{equation}\tag{3.1$c$}\label{3.1c}
  \displaystyle
  v_{h}(x,t) = \frac{1}{h} \int_{t}^{t+h} \widetilde{v}(x,s) \, ds, \ \ \forall \ x \in E \! \setminus \! {Z}_{*}, \ \ \forall \ t \in I.
\end{equation}
\end{lem}

\vspace{0.3cm} \noindent
{\bf Proof.}
We start with the proof of (\ref{3.1a}).

\vspace{0.3cm} \noindent
{\it Case I:} $1 \leq r < \infty$.

\vspace{0.3cm}
Let $E_N \to E$, with $|E_N| < \infty \ \forall \ N$ (if $|E|< \infty$, then simply take $E_N = E$ ahead).
Then, If $1 \leq q < \infty$, we have
\begin{equation*}
\begin{aligned}
 \displaystyle
 \int_{E_{N}} \biggl[ \int_{J} |v(x,t)| dt \biggr] dx
 & = \int_{J} \biggl[ \int_{E_{N}} |v(x,t)| dx \biggr] dt \\
 & \leq \int_{J} {\biggl( \int_{E_{N}} {|v(x,t)|}^{q} dx \biggr)}^{\frac{1}{q}} . {|E_N|}^{1 - \frac{1}{q}} dt \\
 & \leq {|E_N|}^{1 - \frac{1}{q}} . {\Biggl(\int_{J} {\biggl( \int_{E_{N}} {|v(x,t)|}^{q} dx \biggr)}^{\frac{r}{q}} dt\Biggr)}^{\frac{1}{r}} . {|J|}^{1-\frac{1}{r}} \\
 & = {|E_N|}^{1 - \frac{1}{q}} . {|J|}^{1-\frac{1}{r}} {\biggl( \int_{J} {\|v(\cdot,t)\|}^{r}_{L^{q}{(E)}} dt \biggr)}^{\frac{1}{r}} < \infty,
\end{aligned}
\end{equation*}
i.e., (by Fubini and H\"{o}lder) we have
\begin{equation*}
\begin{aligned}
 \displaystyle
 \int_{E_{N}} \biggl[ \int_{J} |v(x,t)| \, dt \biggr] dx < \infty.
\end{aligned}
\end{equation*}
This gives that there exists $Z_N \subseteq E_N$, with $|Z_N| = 0$, such that
\begin{equation*}
 \displaystyle
 \int_{J} |v(x,t)| \, dt < \infty, \ \ \forall \ x \in E_N \! \setminus \! Z_N.
\end{equation*}
In particular, setting
$\displaystyle Z := \bigcup_{N=1}^{\infty} Z_N$,
we have $Z \subseteq E$, $|Z| = 0$, and $\displaystyle \int_{J} |v(x,t)| \, dt < \infty$, for all $x \in E \! \setminus \! Z$.
This shows (\ref{3.1a}) for $1 \leq r < \infty$ and $1 \leq q < \infty$.

\vspace{0.3cm}
If $q = \infty$, we may proceed as follows: taking again $E_N \to E$, with $|E_N| < \infty, \ \forall \ N$, we have:

\begin{equation*}
\begin{aligned}
 \displaystyle
 \int_{E_{N}} \biggl[ \int_{J} |v(x,t)| dt \biggr] dx
 & = \int_{J} \biggl[ \int_{E_{N}} |v(x,t)| dx \biggr] dt \\
 & \leq \int_{J} \biggl( \int_{E_{N}} {\|v(\cdot,t)\|}_{L^{\infty}(E)} dx \biggr) dt \\
 & = \int_{J} {\|v(\cdot,t)\|}_{L^{\infty}(E)} . |E_N| dt \\
 & \leq |E_N| {\biggl(\int_{J} {\|v(\cdot,t)\|}^{r}_{L^{\infty}(E)} dt \biggr)}^{\frac{1}{r}} . {|J|}^{1-\frac{1}{r}} \\
 & \leq |E_N| . {|J|}^{1-\frac{1}{r}} {\biggl(\int_{J} {\|v(\cdot,t)\|}^{r}_{L^{\infty}(E)} dt \biggr)}^{\frac{1}{r}} < \infty,
\end{aligned}
\end{equation*}
i.e.,
\begin{equation*}
\begin{aligned}
 \displaystyle
 \int_{E_{N}} \biggl[ \int_{J} |v(x,t)| dt \biggr] dx < \infty.
\end{aligned}
\end{equation*}
This gives that there exists $Z_N \subseteq E_N$, with $|Z_N| = 0$, such that
$\displaystyle \int_{J} |v(x,t)| \, dt < \infty$, $\ \forall \ x \in E_N \! \setminus \!  Z_N$.
As before, setting $\displaystyle Z := \bigcup_{N=1}^{\infty} Z_N $, it follows that $Z \subseteq E$, $|Z| = 0$, and
$\displaystyle \int_{J} |v(x,t)| \, dt < \infty$, $\ \forall \ x \in E \! \setminus \! Z$, which shows (\ref{3.1a}) for $1 \leq r < \infty$ and $q = \infty$.

\vspace{0.5cm} \noindent
{\it Case II:} $r = \infty$.

As before, consider $E_N \to E$, with $|E_N| < \infty$ $\ \forall \ N$ (if $|E|< \infty$, then simply take $E_N = E$).
Then, if $1 \leq q < \infty$, we have
\begin{equation*}
\begin{aligned}
 \displaystyle
 \int_{E_{N}} \biggl[ \int_{J} |v(x,t)| dt \biggr] \, dx
 & = \int_{J} \biggl[ \int_{E_{N}} |v(x,t)| dx \biggr] \, dt \\
 & \leq \int_{J} {\biggl( \int_{E_{N}} {|v(x,t)|}^{q} \, dx \biggr)}^{\frac{1}{q}} . {|E_N|}^{1 - \frac{1}{q}} \, dt \\
 & \leq {|E_N|}^{1 - \frac{1}{q}} . \int_{J} {\|v(\cdot,t)\|}_{L^{q}(E)} \, dt \\
 & \leq {|E_N|}^{1 - \frac{1}{q}} \int_{J} M \, dt = M . |J| . {|E_N|}^{1 - \frac{1}{q}} < \infty,
\end{aligned}
\end{equation*}
where $\displaystyle M = {\sup{\mathrm{ess}}}_{t \in I} {\| v(\cdot,t) \|}_{L^{q}(E)}$.
Therefore,
\begin{equation*}
 \displaystyle
 \int_{J} |v(x,t)| \, dt \in L^{1}(E_N)
\end{equation*}
(i.e., $\displaystyle \int_{E_N} \bigl( \int_{J} |v(x,t)| \, dt \bigr) dx < \infty$),
and so there must exists $Z_N \subseteq E_N$, with $|Z_N| = 0$, such that
\begin{equation*}
 \displaystyle
 \int_{J} |v(x,t)| \, dt < \infty, \ \ \forall \ x \in E_N \! \setminus \! Z_N.
\end{equation*}
Setting $\displaystyle Z := \bigcup_{N=1}^{\infty} Z_N $, then we have $Z \subseteq E$, $|Z| = 0$ and $\displaystyle \int_{J} |v(x,t)| \, dt < \infty$, $\ \forall \ x \in E \! \setminus \! Z$. This shows (\ref{3.1a}) for $r = \infty$ and $1 \leq q < \infty$.

\vspace{0.3cm}
Now, consider the remaining case $q = \infty$. Taking (again) $E_N \to E$, with $|E_N| < \infty$, $\ \forall \ N$, or simply $E_N = E$, if $|E| < \infty$, setting $\displaystyle M = {\sup{\mathrm{ess}}}_{t \in I} {\| v(\cdot,t) \|}_{L^{\infty}(E)}$, we have that
\begin{equation*}
\begin{aligned}
 \displaystyle
 \int_{E_{N}} \biggl[ \int_{J} |v(x,t)| \, dt \biggr] \, dx
 & = \int_{J} \biggl[ \int_{E_{N}} |v(x,t)| \, dx \biggr] \, dt \\
 & \leq \int_{J} \biggl( \int_{E_{N}} {\|v(\cdot,t)\|}_{L^{\infty}(E_N)} \, dx \biggr) \, dt \\
 & \leq |E_N| \int_{J} {\|v(\cdot,t)\|}_{L^{\infty}(E)} \, dt \\
 & \leq |E_N| \int_{J} M \, dt = M . |E_N| . |J| < \infty,
\end{aligned}
\end{equation*}
i.e., we have
\begin{equation*}
\begin{aligned}
 \displaystyle
 \int_{E_{N}} \biggl[ \int_{J} |v(x,t)| \, dt \biggr] \, dx < \infty.
\end{aligned}
\end{equation*}
Therefore, there exists $Z_N \subseteq E_N$, with $|Z_N| = 0$, such that
$\displaystyle \int_{J} |v(x,t)| \, dt < \infty$, for all $x \in E_N \! \setminus \! Z_N$.
Setting $\displaystyle Z := \bigcup_{N=1}^{\infty} Z_N$, we then have $Z \subseteq E$, $|Z| = 0$, and
\begin{equation*}
 \displaystyle
 \int_{J} |v(x,t)| \, dt < \infty, \ \ \forall \ x \in E \! \setminus \! Z.
\end{equation*}
This shows (\ref{3.1a}) when $r = \infty$ and $q = \infty$; and completes the proof of (\ref{3.1a}).

\vspace{0.3cm}
Proof of (\ref{3.1b}):

Let $v(\cdot,t) \in L^{r}\bigl(I,L^{q}(E)\bigr)$, where $1 \leq q,r \leq \infty$, and $E \in \M({\R}^{n})$ for some interval $I \subseteq \R$. Let $\widetilde{I} = \R$ and
\begin{equation*}
\displaystyle
 \widetilde{v}(\cdot,t) =
 \left\{
  \begin{array}{ccl}
   v(\cdot,t), & \text{se} & t \in I; \\
   0, & \text{se} & t \in \widetilde{I} \backslash I,
  \end{array}
 \right.
\end{equation*}
we then have
\begin{equation*}
 \displaystyle
 \widetilde{v}(\cdot,t) \ \in \ L^{r}\bigl(\widetilde{I},L^{q}(E)\bigr).
\end{equation*}

Taking ${\widetilde{I}}_{l} = \bigl[ {\widetilde{a}}_{l},{\widetilde{b}}_{l} \bigr]$ and making ${\widetilde{I}}_{l} \to \R$ as $l \to \infty$ (i.e., as ${\widetilde{a}}_{l} \to - \infty$ and ${\widetilde{b}}_{l} \to + \infty$), by (\ref{3.1a}) we have that there exists ${Z}_{l} \subseteq E$, with $|{Z}_{l}|=0$, such that
\begin{equation*}
 \displaystyle
 \int_{{\widetilde{I}}_{l}} |\widetilde{v}(x,t)| \, dt < \infty, \ \ \forall \ x \in E \! \setminus \! {Z}_{l}.
\end{equation*}
Taking $\displaystyle {Z}_{*} := \bigcup^{\infty}_{l = 1} {Z}_{l}$, we then have ${Z}_{*} \subseteq E$, with $|{Z}_{*}| = 0$, and
\begin{equation}\tag{3.1$d$}\label{3.1d}
 \displaystyle
 \int_{{\widetilde{I}}_{l}} |\widetilde{v}(x,t)| \, dt < \infty, \ \ \forall \ x \in E \! \setminus \! {Z}_{*}, \ \ \forall \ l.
\end{equation}

Now, given $\widehat{t} \in I$ and $h > 0$ arbitrary, taking $\widehat{l} \in \N$ large enough so that $\bigl[ \widehat{t},\widehat{t}+h \bigr] \subseteq {\widetilde{I}_{\widehat{l}}}$,
we then get for every $x \in E \! \setminus \! {Z}_{*}$, by (\ref{3.1d}), that
\begin{equation*}
 \displaystyle
 \int_{\widehat{t}}^{\widehat{t}+h} |\widetilde{v}(x,t)| \, dt \leq
 \int_{\widetilde{I}_{\, \widehat{l}}} |\widetilde{v}(x,t)| \, dt < \infty.
\end{equation*}
This shows that
\begin{equation*}
 \displaystyle
 \int_{t}^{t+h} |\widetilde{v}(x,t)| \, ds < \infty, \ \ \forall \ x \in E \! \setminus \! {{Z}_{*}}, \ \ \forall \ t \in I, \ \ \forall \ h > 0
\end{equation*}
(with ${Z}_{*} \subseteq E$; $|{Z}_{*}|=0$, and with ${Z}_{*}$ independent of $t \in I$ and of $h>0$). This completes the proof of (\ref{3.1b}).

\vspace{0.3cm}
As an immediate consequence, we obtain the validity of (\ref{3.1c}).

\vspace{0.2cm}
\vspace{-0.75cm}
\begin{flushright}
(Lemma \ref{3.1}) \ $\square$
\end{flushright}

\vspace{0.3cm}
Observe, in the proof of (\ref{3.1a}), that we also have proved the following:

If $v(\cdot,t) \in L^{r}\bigl(I,L^{q}(E)\bigr)$, for some $1 \leq q,r \leq \infty$ (where $E \in \M({\R}^{n})$ and $I \subseteq \R$ is a interval), for each $J \subseteq I$ bounded and each ${E}_{N} \subseteq E$, with $|{E}_{N}| < \infty$, we have
\begin{equation}\tag{3.1$e$}\label{3.1e}
 \displaystyle
 \int_{{E}_{N}} \biggl( \int_{J} |v(x,t)| \, dt \biggr) \, dx < \infty.
\end{equation}

\vspace{0.2cm}
From (\ref{3.1e}), it follows that there exists ${Z}_{t} \subseteq I$ and ${Z}_{*} \subseteq E$, with $|{Z}_{t}|=0$ and $|{Z}_{*}|=0$, such that
\begin{equation}\tag{3.1$f$}\label{3.1f}
 \displaystyle
 \int_{J} |v(x,t)| \, dt < \infty, \ \ \forall \ x \in E \! \setminus \! {Z}_{*}, \text{ and } \forall \ J \subseteq I, \text{ with } J \text{ bounded,}
\end{equation}
and
\begin{equation}\tag{3.1$g$}\label{3.1g}
 \displaystyle
 \int_{K} |v(x,t)| \, dx < \infty, \ \ \forall \ t \in I \! \setminus \! {Z}_{t}, \ \ \forall \ K \subseteq E, \text{ with } |K| < \infty.
\end{equation}

%
%

\section{About the differentiability of the Steklov averages}\label{4}


\vspace{0.3cm}
Let us now relate some properties concerning the differentiation of Steklov averages.

For some $1 \leq {q}_0,{r}_0 \leq \infty$, consider now
\begin{equation*}
 \displaystyle
 v(\cdot,t) \in L^{r_0} \bigl( I,L^{q_0}_{loc}(\Omega) \bigr),
\end{equation*}
where $I \subseteq \R$ is an interval and $\Omega \subseteq {\R}^{n}$ is an (arbitrary) open set.

In particular, for each $K \subseteq \Omega$ compact set, by Lemma \ref{3.1} (see also (\ref{3.1e}), (\ref{3.1f}) and (\ref{3.1g})), we have
\begin{equation*}
 \displaystyle
 v(\cdot,t) \in L^{r_0}(I,L^{q_0}(K)).
\end{equation*}
It follows that there exists ${Z}_{*} \subseteq \Omega$, with $|{\Z}_{*}|=0$, such that
\begin{equation}\tag{4$a$}\label{4a}
 \int_{J} |v(x,t)| \, dt < \infty, \ \ \forall \ x \in \Omega \! \setminus \! {Z}_{*},
\end{equation}
for every bounded interval $J \subseteq I$; and there exists ${Z}_{t} \subseteq I$, with $|{Z}_{t}|=0$, such that
\begin{equation}\tag{4$a$\,'}\label{4a'}
 \displaystyle
 \int_{K} |v(x,t)| \, dx < \infty, \ \ \forall \ t \in I \! \setminus \! {Z}_{t},
\end{equation}
for every compact set $K \subseteq \Omega$.

Setting ${v}_{h}(x,t)$ by
\begin{equation*}
 \displaystyle
 {v}_{h}(x,t) = \frac{1}{h} \int_{t}^{t+h} \widetilde{v}(x,t) \, ds, \ \ \forall \ x \in \Omega \! \setminus \! {Z}_{*}, \ \ \forall \ t \in I
\end{equation*}
(where $h > 0$ is given, and $\widetilde{v}(\cdot,t) = v(\cdot,t)$, if $t \in I$; or $\widetilde{v}(\cdot,t)=0$, if $t \in {\R} \backslash I$), we have
\begin{equation}\tag{4$b$}\label{4b}
 \displaystyle
 {v}_{h}(\cdot,t) \in C^{0}\bigl(I,L^{q_0}_{loc}(\Omega)\bigr) \, \cap \, L^{\infty}\bigl(I,L^{q_0}_{loc}(\Omega)\bigr).
\end{equation}
Note that (\ref{4b}) means that for each compact $K \in \Omega$, one has
\begin{equation*}
 \displaystyle
 {v}_{h}(\cdot,t) \in C^{0}\bigl(I,L^{q_0}(K)\bigr) \, \cap \, L^{\infty}\bigl(I,L^{q_0}(K)\bigr).
\end{equation*}

Now, consider $v(\cdot,t) \in L^{r_0}\bigl(I,L^{q_0}_{loc}(\Omega)\bigr)$, with $1 \leq {q}_0, {r}_0 \leq \infty$, such that
\begin{equation}\tag{4$c$}\label{4c}
 \displaystyle
 \nabla v(\cdot,t) \in L^{r_1}\bigl(I,L^{q_1}_{loc}(\Omega)\bigr) \ \ (\text{with } 1 \leq {q}_1, {r}_1 \leq \infty),
\end{equation}
where $\nabla v(\cdot,t)$ is meant in the {\it distributional sense}: for each $t \in I \! \setminus \! {Z}_{t}$ (with $|{Z}_{t}|=0$), by (\ref{4a'}) we have
$v(\cdot,t) \in L^{1}_{loc}(\Omega)$.
In particular, we can compute its distributional derivative ${D}_{i}v(\cdot,t)$, for $1 \leq i \leq n$, which are given by:
\begin{equation*}
 \begin{array}{rcl}
  \Bigl\langle {D}_{i}v(\cdot,t) \, \big| \, \phi \Bigr\rangle
  & = & - \Bigl\langle v(\cdot,t) \, \big| \, \frac{\partial \phi}{\partial {x}_{i}} \Bigr\rangle \\
  & = & \displaystyle - \int_{\Omega} v(x,t) \frac{\partial \phi}{\partial {x}_{i}} (x) \, dx, \ \ \forall \ \phi \in C^{\infty}_{0}(\Omega).
 \end{array}
\end{equation*}

For $1 \leq i \leq n$, the assumption (\ref{4c}) says that for almost all $t \in I$, ${D}_{i}v(\cdot,t)$ is given by some function ${g}_{i}(\cdot,t) \in L^{q_1}_{loc}(\Omega)$, and we have
\begin{equation*}
 \left\{
  \begin{array}{lcl}
   \displaystyle \int_{I} {\bigl\| {g}_{i}(\cdot,t) \bigr\|}^{r_1}_{L^{q_1}(K)} dt < \infty, & \text{se} & 1 \leq {r}_{1} < \infty, \\
   \displaystyle {\sup{\mathrm{ess}}}_{t \in I} \, {\|{g}_{i}(\cdot,t)\|}_{L^{q_1}(K)} < \infty, & \text{se} & {r}_{1} = \infty,
  \end{array}
 \right.
\end{equation*}
for each given compact set $K \subseteq \Omega$.
Follows then from (\ref{4c}) (by Lemma \ref{2.4} part (\ref{2.4b})) that we have
\begin{equation*}
 \displaystyle
 {({D}_{i}v)}_{h}(\cdot,t) \in C^{0}\bigl(I,L^{q_1}_{loc}(\Omega)\bigr) \, \cap \, L^{\infty}\bigl(I,L^{q_1}_{loc}(\Omega)\bigr), \ \ 1 \leq i \leq n.
\end{equation*}

Moreover, by Lemma \ref{3.1}, we have (enlarging the null sets ${Z}_{t}$ and ${Z}_{*}$, if necessary) that there exists ${Z}_{t} \subseteq I$, with $|{Z}_{t}|=0$, such that
\begin{equation*}
 \displaystyle \int_{K} |v(x,t)| \, dx < \infty, \ \ \forall \ t \in I \! \setminus \! {Z}_{t}, \text{ for each compact set } K \subseteq \Omega,
\end{equation*}
and
\begin{equation*}
 \displaystyle \int_{K} \big| {D}_{i}v(x,t) \big| \, dx < \infty, \ \ \forall \ t \in I \! \setminus \! {Z}_{t}, \ \ \forall \ 1 \leq i \leq n, \ \ \forall \ K \subseteq \Omega;
\end{equation*}
and there exists ${Z}_{*} \subseteq \Omega$, with $|{Z}_{*}|=0$,
such that
\begin{equation*}
 \displaystyle
 \int_{I} |v(x,t)| \, dt < \infty, \ \ \forall \ x \in \Omega \! \setminus \! {Z}_{*},
\end{equation*}
and
\begin{equation*}
 \displaystyle
 \int_{I} \big| {D}_{i}v(x,t) \big| \, dt < \infty, \ \ \forall \ x \in \Omega \! \setminus \! {Z}_{*}, \ \text{ (}\forall \ 1 \leq i \leq n \text{)}.
\end{equation*}

In particular, for the pointwise values of ${v}_{h}(\cdot,t)$ and of ${\bigl({D}_{i}v\bigr)}_{h}(\cdot,t)$ we have:
\begin{equation*}
 \displaystyle
 {v}_{h}(x,t) = \frac{1}{h} \int_{t}^{t+h} \widetilde{v}(x,s) \, ds, \ \ \forall \ x \in \Omega \! \setminus \! {Z}_{*}, \ \ \forall \ t \in I;
\end{equation*}
\begin{equation*}
 \begin{array}{rcl}
  {\bigl({D}_{i}v\bigr)}_{h}(x,t)
  & = & \displaystyle \frac{1}{h} \int_{t}^{t+h} \widetilde{{D}_{i}v}(x,s) \ ds, \ \ \forall \ x \in \Omega \! \setminus \! {Z}_{*}, \ \ \forall \ t \in I \\
  \Biggl(
  & = & \displaystyle \frac{1}{h} \int_{t}^{t+h} \widetilde{g_{i}}(x,s) \, ds, \, \text{ where }
   \widetilde{g_{i}}(x,s) =
   \left\{
    \begin{array}{l}
     {g}_{i}(\cdot,t), \ t \in I, \\
     0, \ t \in {\R} \! \setminus \! I.
   \end{array}
  \right.
  \Biggr)
 \end{array}
\end{equation*}
and, more importantly,
\begin{equation*}
\displaystyle
 {\bigl({D}_{i}v\bigr)}_{h}(x,t) = \displaystyle \frac{1}{h} \int_{t}^{t+h} \! \bigl({D}_{i}v\bigr)(x,s) \, ds, \quad \forall \ x \in \Omega \! \setminus \! {Z}_{*}, \ \ \forall \ t \in {I}_{h};
\end{equation*}
and, for each compact set $K \subseteq \Omega$, with:
\begin{equation*}
 \displaystyle
 {v}_{h}(\cdot,t) \ \in \ C^{0}\bigl(I,L^{q_0}(K)\bigr) \, \cap \, L^{\infty}\bigl(I,L^{q_0}(K)\bigr); \ \text{ and }
\end{equation*}
\begin{equation*}
 \displaystyle
 {\bigl({D}_{i}v\bigr)}_{h}(\cdot,t) \ \in \ C^{0}\bigl(I,L^{q_1}(K)\bigr) \, \cap \, L^{\infty}\bigl(I,L^{q_1}(K)\bigr).
\end{equation*}

%
%

\vspace{0.3cm}
In the next result we'll shows that the operators ${D}_{i}$ and ${(\, \cdot \,)}_{h}$ commute.

\vspace{0.3cm}
\begin{lem}\label{4.1}
Let $I \subseteq \R$ an interval and $\Omega \subseteq {\R}^{n}$ an open set.
Let $ v(\cdot,t) \in L^{r_0}\bigl(I,L^{q_0}_{loc}(\Omega)\bigr)$ (for some $1 \leq {q}_{0},{r}_{0} \leq \infty)$, such that $\nabla v(\cdot,t) \in L^{r_1}\bigl(I,L^{q_1}_{loc}(\Omega)\bigr)$ (for some $1 \leq {q}_{1},{r}_{1} \leq \infty)$.

Then, for each $1 \leq i \leq n$ we have
\begin{equation}\tag{4.1}\label{4.1}
 \displaystyle
 {D}_{i}\bigl({v}_{h}(\cdot,t)\bigr) = {\bigl({D}_{i}v\bigr)}_{h}(\cdot,t), \ \ \forall \ t \in \ {I}_{h}, \text{ where } {I}_{h} = \bigl\{ t \in I \ | \ t+h \in \overline{I} \bigr\}.
\end{equation}
\end{lem}

\vspace{0.3cm}
Observe that: ${v}_{h} \in L^{q_0}_{loc}(\Omega), \ \forall \ t \in I$;
and ${\bigl({D}_{i}v\bigr)}_{h}(\cdot,t) \in L^{q_1}_{loc}(\Omega), \ \forall \ t \in I$. In particular, under the hypothesis of Lemma \ref{4.1}, we have
\begin{equation*}
 \displaystyle
 {v}_{h}(\cdot,t) \in C^{0}\bigl(I,L^{q_0}(K)\bigr) \, \cap \, L^{\infty}\bigl(I,L^{q_0}(K)\bigr), \ \text{ and }
\end{equation*}
\begin{equation*}
 \displaystyle
 \nabla {v}_{h}(\cdot,t) \in C^{0}\bigl({I}_{h},L^{q_1}(K)\bigr) \, \cap \, L^{\infty}\bigl({I}_{h},L^{q_1}(K)\bigr)
\end{equation*}
where ${I}_{h} = \{ t \in I \ | \ t+h \in \overline{I} \}$, for each compact ${K} \in \Omega$.

\vspace{0.3cm} \noindent
{\bf Proof of Lemma \ref{4.1}.}
Given $t \in I$, with $t+h \in \overline{I}$, and $\phi \in C^{\infty}_{0}(\Omega)$ such that $supp(\phi) \subseteq K \subseteq \Omega$, where $K$ is an compact, we have ${v}_{h}(\cdot,t) \in L^{q_0}_{loc}(\Omega) \subseteq {\cal{D}}'(\Omega)$
and
\begin{equation*}
\begin{aligned}
 \displaystyle
 \Bigl\langle {D}_{i} \bigl( {v}_{h}(\cdot,t) \bigr) \ \big| \ \phi \Bigr\rangle
 & = - \Bigl\langle {v}_{h}(\cdot,t) \ \big| \ {D}_{i} \phi \Bigr\rangle
 =
 - \int_{K} {v}_{h}(x,t) \frac{\partial \phi}{\partial {x}_{i}}(x) \, dx \\
 & = - \frac{1}{h} \! \int_{K} \! \biggl( \int_{t}^{t+h} v(x,s) \frac{\partial \phi}{\partial {x}_{i}}(x) ds \biggr) dx
 =
 - \frac{1}{h} \! \int_{t}^{t+h} \! \biggl( \int_{\Omega} v(x,s) \frac{\partial \phi}{\partial {x}_{i}}(x) dx \biggr) ds \\
 & = - \frac{1}{h} \int_{t}^{t+h} \biggl\langle v(\cdot,s) \ \big| \ \frac{\partial \phi}{\partial {x}_{i}} \biggr\rangle \, ds
 =
 \frac{1}{h} \int_{t}^{t+h} \Bigl\langle {D}_{i}v(\cdot,s) \ \big| \ \phi \Bigr\rangle \, ds \\
 & = \frac{1}{h} \int_{t}^{t+h} \Bigl\langle {g}_{i}(\cdot,s) \ \big| \ \phi \Bigr\rangle \, ds
 =
 \frac{1}{h} \int_{t}^{t+h} \biggl( \int_{\Omega} {g}_{i}(x,s) \ \phi(x) \, dx \biggr) \, ds \\
 & = \frac{1}{h} \int_{K} \phi(x) \, \biggl( \int_{t}^{t+h} {g}_{i}(x,s) \, ds \biggr) \, dx
 =
 \frac{1}{h} \int_{K} \phi(x) \, \biggl( \int_{t}^{t+h} \bigl({D}_{i}v\bigr)(x,s) \, ds \biggr) \, dx \\
 & = \int_{K} {\bigl({D}_{i}v\bigr)}_{h}(x,t) \cdot \phi(x) \, dx
 =
 \Bigl\langle {\bigl({D}_{i}v\bigr)}_{h}(\cdot,t) \ \big| \ \phi \Bigr\rangle,
\end{aligned}
\end{equation*}
i.e.,
\begin{equation*}
 \displaystyle
 \Bigl\langle {D}_{i} \bigl( {v}_{h}(\cdot,t) \bigr) \ \big| \ \phi \Bigr\rangle = \Bigl\langle {\bigl({D}_{i}v\bigr)}_{h}(\cdot,t) \ \big| \ \phi \Bigr\rangle, \ \  \forall \ t \in I, \text{ with } t+h \in \overline{I},
\end{equation*}
where $\phi \in C^{\infty}_{0}(\Omega)$ is an arbitrary test function, and $\langle \cdot,\cdot \rangle$ denotes the natural pairing of ${\cal{D}}'(\Omega)$ and ${\cal{D}}(\Omega)$.
This shows that
\begin{equation*}
 \displaystyle
 {D}_{i}\bigl( {v}_{h}(\cdot,t) \bigr) = {\bigl( {D}_{i}v \bigr)}_{h}(\cdot,t), \ \ t \in {I}_{h} = \{ t \in I \ | \ t+h \in \overline{I} \}.
\end{equation*}

\vspace{-0.5cm}
\begin{flushright}
(Lemma \ref{4.1}) \ $\square$
\end{flushright}

%
%

\vspace{0.3cm}
Another important operation is the differentiability of ${v}_{h}(\cdot,t) \in L^{q}(E)$ with respect to $t$ (in the Banach space $L^{q}(E)$):

\vspace{0.3cm}
\begin{lem}\label{4.2}
Given an interval $I \subseteq \R$, $E \in {\cal{M}}({\R}^{n})$, and $1 \leq q \leq \infty$; let
\begin{equation*}
 \displaystyle
 v(\cdot,t) \in C^{0}\bigl(I,L^{q}(E)\bigr).
\end{equation*}
Let ${I}_{h} = \bigl\{ t \in I \ | \ (t+h) \in I \bigr\}$ and ${\ring{I}}_{h} = int({I}_{h}) = \bigl\{ t \in \ring{I} \ | \ (t+h) \in \ring{I} \bigr\}$,
where $h>0$ (small enough that ${\ring{I}}_{h} \neq \emptyset$) is given.

Then, for every $t \in {\ring{I}}_{h}$, we have
\begin{equation}\tag{4.2$a$}\label{4.2a}
 \displaystyle
  {\biggl\| \frac{{v}_{h}(\cdot,t+\Delta t) - {v}_{h}(\cdot,t)}{\Delta t} - \frac{v(\cdot,t+h)-v(\cdot,t)}{h} \biggr\|}_{L^{q}(E)} \longrightarrow 0,
  \ \text{ as } \Delta t \to 0,
\end{equation}
uniformly on $t \in [t_1,t_2] \subseteq {\ring{I}}_{h}$.
\end{lem}

\vspace{0.3cm}
In other words:
\begin{equation*}
 \displaystyle
 {v}_{h}(\cdot,t) = \frac{1}{h} \int_{t}^{t+h} v(\cdot,s) \, ds \ \ \ (\text{for }t \in {I}_{h})
\end{equation*}
is (pointwise) strongly differentiable at $t \in {\ring{I}}_{h}$, with
\begin{equation*}
 \displaystyle
 {\bigl({v}_{h}\bigr)}_{t}(\cdot,t) = \frac{1}{h}\biggl( v(\cdot,t+h) - v(\cdot,t) \biggr) \ \ \ \bigl(\in L^{q}(E)\bigr), \ \forall \ t \in {\ring{I}}_{h}.
\end{equation*}

\vspace{0.3cm} \noindent
{\bf Proof of Lemma \ref{4.2}.}
Given $\widehat{t} \in [{t}_{1},{t}_{2}]$, with $\bigl[{t}_{1},{t}_{2}\bigr] \subseteq {\ring{I}}_{h}$ an compact interval, take ${\delta}_{0} > 0$ small enough that
\begin{equation*}
 \displaystyle
 \bigl[{t}_{1} - {\delta}_{0} \, , \, {t}_{2} + h + {\delta}_{0}\bigr] \ \subseteq \ {\ring{I}}_{h}.
\end{equation*}

Given $\varepsilon > 0$, $\bigl($because $v(\cdot,t)$ is uniformly continuous on $[{t}_{1} - {\delta}_{0} \, , \, {t}_{2} + h + {\delta}_{0}]\bigr)$ we can take $\delta = \delta(\varepsilon) \leq {\delta}_{0}$ small enough that
\begin{equation*}
 \displaystyle
 {\bigl\| v(\cdot,t) - v(\cdot,s) \bigr\|}_{L^{q}(E)} \leq \varepsilon, \ \ \forall \ s,t \in [{t}_{1} \! - \! {\delta}_{0} , {t}_{2} \! + \! h + \! {\delta}_{0}], \text{ with } |s-t|\leq \delta.
\end{equation*}

Then, for any $\widehat{t} \in [{t}_{1},{t}_{2}]$ and any $\Delta t \in \R$ with $0 < |\Delta t| < \delta$, we have
\begin{equation*}
\begin{aligned}
 \displaystyle
 & {\Bigl\| \frac{1}{\Delta t} \Bigl( {v}_{h}(\cdot,\widehat{t} + \Delta t) - {v}_{h}(\cdot, \widehat{t}) \Bigr) - \frac{1}{h} \Bigl( v(\cdot,\widehat{t}+h) - v(\cdot,\widehat{t}) \Bigr) \Bigr\|}_{L^{q}(E)} = \\
 & \quad = \frac{1}{h} \ {\Bigl\| \frac{1}{\Delta t} \int_{\widehat{t} + \Delta t}^{\widehat{t} + \Delta t + h} \hspace{-0.5cm} v(\cdot,s) \, ds - \frac{1}{\Delta t} \int_{\widehat{t}}^{\widehat{t} + h} \hspace{-0.5cm} v(\cdot,s) \, ds \ - \ v(\cdot,\widehat{t} + h) \ + \ v(\cdot,\widehat{t}) \Bigr\|}_{L^{q}(E)} \\
 & \quad = \frac{1}{h} \ {\Bigl\| \Bigl( \frac{1}{\Delta t} \int_{\widehat{t} + h}^{\widehat{t} + \Delta t + h} \hspace{-0.5cm} v(\cdot,s) \, ds - v(\cdot,\widehat{t} + h) \Bigr) - \Bigl( \frac{1}{\Delta t} \int_{\widehat{t}}^{\widehat{t} + \Delta t} \hspace{-0.5cm} v(\cdot,s) \, ds \ - \ v(\cdot,\widehat{t}) \Bigr) \Bigr\|}_{L^{q}(E)} \\
 & \quad \leq \frac{1}{h} \ {\Bigl\| \frac{1}{\Delta t} \int_{\widehat{t} + h}^{\widehat{t} + \Delta t + h} \hspace{-0.5cm} v(\cdot,s) \, ds - v(\cdot,\widehat{t} + h) \Bigr\|}_{L^{q}(E)} \! \! \! \! + \frac{1}{h} \ {\Bigl\| \frac{1}{\Delta t} \int_{\widehat{t}}^{\widehat{t} + \Delta t} \hspace{-0.5cm} v(\cdot,s) \, ds \ - \ v(\cdot,\widehat{t}) \Bigr\|}_{L^{q}(E)} \\
 & \quad = \frac{1}{h} \ {\Bigl\| \frac{1}{\Delta t} \int_{\widehat{t} + h}^{\widehat{t} + \Delta t + h} \! \! \Bigl( v(\cdot,s) - v(\cdot,\widehat{t}+h) \Bigr) ds \Bigr\|}_{L^{q}(E)} \! \! \! \! + \frac{1}{h} \ {\Bigl\| \frac{1}{\Delta t} \int_{\widehat{t}}^{\widehat{t} + \Delta t} \! \Bigl( v(\cdot,s) - v(\cdot,\widehat{t}) \Bigr) ds \Bigr\|}_{L^{q}(E)} \\
 & \quad \leq \frac{1}{h} \biggl( \frac{1}{\Delta t} \int_{\widehat{t} + h}^{\widehat{t} + \Delta t + h} \! \! {\bigl\| v(\cdot,s) - v(\cdot,\widehat{t}+h) \bigr\|}_{L^{q}(E)} \! ds \biggr) \! + \frac{1}{h} \biggl( \frac{1}{\Delta t} \int_{\widehat{t}}^{\widehat{t} + \Delta t} \! {\bigl\| v(\cdot,s) - v(\cdot,\widehat{t}) \bigr\|}_{L^{q}(E)} \! ds \biggr) \\
 & \quad \leq \frac{1}{h} \biggl( \frac{1}{\Delta t} \int_{\widehat{t} + h}^{\widehat{t} + \Delta t + h} \! \! \! \! \varepsilon \ ds \biggr) + \frac{1}{h} \biggl( \frac{1}{\Delta t} \int_{\widehat{t}}^{\widehat{t} + \Delta t} \! \! \! \! \varepsilon \ ds \biggr) \quad = \quad \frac{1}{h} \ \varepsilon + \frac{1}{h} \ \varepsilon  \quad = \quad \frac{2}{h} \ \varepsilon.
\end{aligned}
\end{equation*}
i.e., for any $\Delta t \in \R$ with $0 < |\Delta t| < \delta$ we have
\begin{equation*}
 \displaystyle
 {\Bigl\| \frac{1}{\Delta t} \Bigl( {v}_{h}(\cdot,\widehat{t} + \Delta t) - {v}_{h}(\cdot, \widehat{t}) \Bigr) - \frac{1}{h} \Bigl( v(\cdot,\widehat{t}+h) - v(\cdot,\widehat{t}) \Bigr) \Bigr\|}_{L^{q}(E)} \leq \frac{2}{h} \ \varepsilon,
\end{equation*}
for all $\widehat{t} \in [{t}_{1},{t}_{2}]$, where $\delta > 0$ depends only on $\varepsilon$ (and not on $\widehat{t} \in \bigl[{t}_{1},{t}_{2}\bigr]$).
This shows (\ref{4.2a}), so that the Lemma \ref{4.2} is now proven.

\vspace{-0.6cm}
\begin{flushright}
(Lemma \ref{4.2}) \ $\square$
\end{flushright}

\vspace{0.3cm}
Hence, we have that the mapping
\begin{equation*}
 \displaystyle
 {v}_{h} : {I}_{h} \to L^{q}(E) \quad \Bigl(\in C^{0}\bigl({I}_{h},L^{q}(E)\bigr)\Bigr)
\end{equation*}
is {\it strongly differentiable} at each $t \in {\ring{I}}_{h}$ (when $v(\cdot,t) \in C^{0}\bigl(I,L^{q}(E)\bigr)$, for $1 \leq q \leq \infty$), with:
\begin{equation}\tag{4.2$b$}\label{4.2b}
 \displaystyle
 {\bigl({v}_{h}\bigr)}_{t}(\cdot,t) \ \equiv \ \frac{\partial}{\partial t} \ {v}_{h}(\cdot,t) \ = \ \frac{v(\cdot,t+h) - v(\cdot,t)}{h}, \quad \forall \ t \in {\ring{I}}_{h}.
\end{equation}

%
%

\vspace{0.3cm}
The next result shows that, in the case where we only have
$v(\cdot,t) \in L^{r}\bigl(I,L^{q}(E)\bigr)$, with $1 \leq r < \infty$,
then (\ref{4.2b}) still holds, but only {\it almost everywhere} on ${\ring{I}}_{h}$.

\vspace{0.3cm}
\begin{lem}\label{4.3}
Given an interval $I \subseteq \R$, $E \in {\cal{M}}({\R}^{n})$, $1 \leq q \leq \infty$, $1 \leq r < \infty$ and $h>0$ (such that ${\ring{I}}_{h}$ is not empty), let
\begin{equation*}
  \displaystyle
  v(\cdot,t) \in L^{r}\bigl(I,L^{q}(E)\bigr)
\end{equation*}
(and, as consequence, ${v}_{h}(\cdot,t) \in C^{0}\bigl({\overline{I}}_{h},L^{q}(E)\bigr)$).

Let ${I}_{h} = \{ t \in I \ | \ (t+h) \in I \}$ and
${\ring{I}}_{h} \ = \ int({I}_{h})$ (as in Lemma \ref{4.2}).

Then, there exists ${Z}_{**} \subseteq I$, with $|{Z}_{**}|=0$, such that
\begin{equation*}
 \displaystyle
 {v}_{h}(\cdot,t) : {\overline{I}}_{h} \to L^{q}(E)
\end{equation*}
is (strongly) differentiable at every $t \in {\ring{I}}_{h} \! \setminus \! {Z}_{**}$, with
\begin{equation}\tag{4.3$a$}\label{4.3a}
 \displaystyle
 {\bigl({v}_{h}\bigr)}_{t}(\cdot,t) = \frac{v(\cdot,t+h) - v(\cdot,t)}{h} \qquad \forall \ t \in {\ring{I}}_{h} \! \setminus \! {Z}_{**}.
\end{equation}
\end{lem}

\vspace{0.3cm}
Lemma \ref{4.3} says that
\begin{equation*}
 \displaystyle
 v(\cdot,t+h) \in L^{q}(E) \ \text{ and } \ v(\cdot,t) \in L^{q}(E), \ \ \forall \ t \in {\ring{I}}_{h} \! \setminus \! {Z}_{**};
\end{equation*}
and that
\begin{equation*}
 \displaystyle
 {v}_{h}(\cdot,t) = \frac{1}{h} \int_{t}^{t+h} v(\cdot,s) \, ds \quad \bigl( \in C^{0}\bigl({\overline{I}}_{h},L^{q}(E)\bigr)
\end{equation*}
is strongly differentiable at each $t \in {\ring{I}}_{h} \! \setminus \! {Z}_{**}$, and (\ref{4.3a}) holds, i.e.,
\begin{equation*}
 \displaystyle
 {\Bigl\| \frac{{v}_{h}(\cdot,t + \Delta t) - {v}_{h}(\cdot,t)}{\Delta t} - \frac{ v(\cdot,t+h) - v(\cdot,t)}{h} \Bigr\|}_{L^{q}(E)} \longrightarrow 0,
\end{equation*}
as $\Delta t \to 0$, for each $t \in {\ring{I}}_{h} \! \setminus \! {Z}_{**}$.

\vspace{0.3cm} \noindent
{\bf Proof of Lemma \ref{4.3}.}
Given $v(\cdot,t) \in L^{r}\bigl(I,L^{q}(E)\bigr)$, let ${Z}_{0,0} \subseteq I$ be such that $|{Z}_{0,0}|=0$ and
\begin{equation*}
 \displaystyle
 v(\cdot,t) \in L^{q}(E), \quad \forall \ t \in I \! \setminus \! {Z}_{0,0}.
\end{equation*}
Let ${Z}_{0,h} \equiv \{ t \in {I}_{h} \ | \ (t+h) \in {Z}_{0,0} \} \subseteq {Z}_{0,0} - h$. Then, we have ${Z}_{0,h} \subseteq {I}_{h}$, with $|{Z}_{0,h}|=0$, and
\begin{equation*}
 v(\cdot,t+h) \in L^{q}(E), \qquad \forall \ \ t \in {I}_{h} \backslash {Z}_{0,h}.
\end{equation*}
$\bigl($In particular, we have $v(\cdot,t)$ and $v(\cdot,t+h)$ in $L^{q}(E)$, $\forall \ t \in {I}_{h} \! \setminus \! \bigl({Z}_{0,0} \cup {Z}_{0,h}\bigr)\bigr)$.

Taking ${Z}_{0} : = {Z}_{0,0} \cup {Z}_{0,h}$, we have ${Z}_{0} \subseteq I$, with $|{Z}_{0}|=0$, and
\begin{equation*}
 \displaystyle
 v(\cdot,t), \ v(\cdot,t+h) \in L^{q}(E), \quad \forall \ t \in {I}_{h} \! \setminus \! {Z}_{0}.
\end{equation*}

Now, because $v(\cdot,t) \in L^{r}\bigl(I,L^{q}(E)\bigr)$ for some $1 \leq r < \infty$, we can take a sequence of smooth approximations
${w}_{m}(\cdot,t) \in C^{0}_{c}\bigl({\ring{I}},L^{q}(E)\bigr)$
such that
\begin{equation}\tag{4.3$b$}\label{4.3b}
 \displaystyle
 {\bigl\| {w}_{m} - v \bigr\|}_{L^{r}(I,L^{q}(E))} \longrightarrow 0, \ \ (\text{ as } m \to \infty),
\end{equation}
and (passing to a subsequence, if necessary):
\begin{equation}\tag{4.3$c$}\label{4.3c}
 \displaystyle
 {\bigl\| {w}_{m}(\cdot,t) - v(\cdot,t) \bigr\|}_{L^{q}(E)} \longrightarrow 0, \ \ (\text{ as } m \to \infty), \ \forall \ t \in I \! \setminus \! {Z}_{*,0},
\end{equation}
for some ${Z}_{*,0} \subseteq I$, with $|{Z}_{*,0}|=0$ and $\ {Z}_{*,0} \supseteq {Z}_{0,0}$.

In particular, setting ${Z}_{*,h} : = \bigl\{ t \in {I}_{h} \ | \ t+h \in {Z}_{*,0} \bigr\} \subseteq {Z}_{*,0} - h$, we have
${Z}_{*,h} \subseteq {I}_{h}$, with $|{Z}_{*,h}|=0$, $\ {Z}_{*,h} \supseteq {Z}_{0,h}$, and
\begin{equation}\tag{4.3$c$\,'}\label{4.3c'}
 \displaystyle
 {\bigl\| {w}_{m}(\cdot,t+h) - v(\cdot,t+h) \bigr\|}_{L^{q}(E)} \hspace{-0.5cm} \longrightarrow 0, \quad (\text{ as } m \to \infty), \ \forall \ \ t \in {I}_{h} \! \setminus \! {Z}_{*,h}.
\end{equation}

In particular, letting ${Z}_{*} \subseteq I$ be given by
\begin{equation*}
 \displaystyle
 {Z}_{*} = {Z}_{*,0} \cup {Z}_{*,h},
\end{equation*}
we have ${Z}_{*} \subseteq I$, with $|{Z}_{*}|=0$ and $\ {Z}_{*} \supseteq {Z}_{0}$; and
\begin{equation*}
 \displaystyle
 {\bigl\| {w}_{m}(\cdot,t) - v(\cdot,t) \bigr\|}_{L^{q}(E)} \hspace{-0.5cm} \longrightarrow 0, \qquad {\bigl\| {w}_{m}(\cdot,t+h) - v(\cdot,t+h) \bigr\|}_{L^{q}(E)} \hspace{-0.5cm} \longrightarrow 0,
\end{equation*}
as $m \to \infty$, $\forall \ t \in {I}_{h} \! \setminus \! {Z}_{*}$.

Finally, for each $m = 1,2,3,...$, by the standard Lebesgue's differentiation theorem, we have that there exists some null set
${Z}_{m,0} \subseteq I$, with $|{Z}_{m,0}|=0$ and ${Z}_{m,0} \supseteq {Z}_{0,0}$, such that
\begin{equation}\tag{4.3$d$}\label{4.3d}
 \displaystyle
 \lim_{\Delta t \to 0} \frac{1}{\Delta t} \int_{t}^{t+\Delta t} {\bigl\| {v}(\cdot,s) - {w}_{m}(\cdot,s) \bigr\|}_{L^{q}(E)}^{r} ds =  {\bigl\| {v}(\cdot,t) - {w}_{m}(\cdot,t) \bigr\|}_{L^{q}(E)}^{r}, \ \ \forall \ t \in I \! \setminus \! {Z}_{m,0}.
\end{equation}

Letting
\begin{equation*}
 \displaystyle
 {Z}_{m,h} : = \bigl\{ t \in {I}_{h} \ | \ (t+h) \in {Z}_{m,0} \bigr\} \subseteq -h + {Z}_{m,0},
\end{equation*}
we have
${Z}_{m,h} \subseteq {I}_{h}$, with $|{Z}_{m,h}|=0$ and $\ {Z}_{m,h} \supseteq {Z}_{0,h}$; and (by (\ref{4.3d})):
\begin{equation}\tag{4.3$e$}\label{4.3e}
 \displaystyle
 \lim_{\Delta t \to 0} \! \frac{1}{\Delta t}  \int_{t+h}^{t + h + \Delta t}  \! {\bigl\| {v}(\cdot,s) - {w}_{m}(\cdot,\! s) \bigr\|}_{L^{q}(E)}^{r} \, ds = {\bigl\| {v}(\cdot,t  + h) - {w}_{m}(\cdot,t + h) \bigr\|}_{L^{q}(E)}^{r},
\end{equation}
$\forall \ \ t \in {I}_{h} \! \backslash {Z}_{m,h}$.
Then, if we taking ${Z}_{**} \subseteq I$ and the null set given by
\begin{equation*}
 \displaystyle
 {Z}_{**} : = {Z}_{*,0} \cup {Z}_{*,h} \cup \biggl( \bigcup_{m=1}^{\infty} \bigl( {Z}_{m,0} \cup {Z}_{m,h} \bigr) \biggr),
\end{equation*}
given (any) $t \in {\ring{I}}_{h} \backslash {Z}_{**}$, we will then have:
\begin{equation}\tag{4.3$f$}\label{4.3f}
 \displaystyle
 \lim_{\Delta t \to 0} {\biggl\| \frac{{v}_{h}(\cdot,t + \Delta t) - {v}_{h}(\cdot,t)}{\Delta t} - \frac{v(\cdot,t + h) - v(\cdot,t)}{h} \biggr\|}_{L^{q}(E)} \! \! \! = 0, \ \ \forall \ t \in {I}_{h} \! \setminus \! {Z}_{**},
\end{equation}
which will shows (\ref{4.3a}), concluding the proof of the Lemma \ref{4.3}.

\vspace{0.3cm} \noindent
{\it Claim:} (\ref{4.3f}) is true.

Indeed, given $\widehat{t} \in {\ring{I}}_{h} \! \setminus \! {Z}_{**}$, we may proceed as follows:
given $\varepsilon > 0$, let $w \equiv {w}_{m} \in C^{0}_{c}\bigl({\ring{I}},L^{q}(E)\bigr)$ be some term of the sequence ${\bigl({w}_{m}\bigr)}_{m}$ given in (\ref{4.3b}), (\ref{4.3c}) and (\ref{4.3c'}), such that we have
\begin{equation*}
 \displaystyle
 {\bigl\| w - v \bigr\|}_{L^{r}(I,L^{q}(E))} \leq \varepsilon,
\end{equation*}
and
\begin{equation*}
 \displaystyle
 {\bigl\| w(\cdot,\widehat{t}\,) - v(\cdot,\widehat{t}\,) \bigr\|}_{L^{q}(E)} \leq \varepsilon,
\end{equation*}
and
\begin{equation*}
 \displaystyle
 {\bigl\| w(\cdot,\widehat{t}+h) - v(\cdot,\widehat{t}+h) \bigr\|}_{L^{q}(E)} \leq \varepsilon.
\end{equation*}

This gives
\begin{equation*}
\begin{aligned}
 \displaystyle
 & \limsup_{\Delta t \to 0} {\Bigl\| \frac{1}{\Delta t} \Bigl( {v}_{h}(\cdot,\widehat{t} + \Delta t) - {v}_{h}(\cdot, \widehat{t}) \Bigr) - \frac{1}{h} \Bigl( v(\cdot,\widehat{t}+h) - v(\cdot,\widehat{t}) \Bigr) \Bigr\|}_{L^{q}(E)} = \\
 & \quad = \frac{1}{h} \limsup_{\Delta t \to 0}  \Bigl\| \frac{1}{\Delta t} \int_{\widehat{t} + \Delta t}^{\widehat{t} + \Delta t + h} \hspace{-0.5cm} v(\cdot,s) \, ds - \frac{1}{\Delta t} \int_{\widehat{t}}^{\widehat{t}+h} \hspace{-0.3cm} v(\cdot,s) \, ds \ - \ v(\cdot,\widehat{t}\!+\!h) \ + \ v(\cdot,\widehat{t}) {\Bigr\|}_{L^{q}(E)} \\
 & \quad = \frac{1}{h} \limsup_{\Delta t \to 0} \biggl\| \biggl( \frac{1}{\Delta t} \int_{\widehat{t} + h}^{\widehat{t} + \Delta t + h} \hspace{-0.5cm} v(\cdot,s) \, ds - v(\cdot,\widehat{t} + h) \biggr) - \biggl( \frac{1}{\Delta t} \int_{\widehat{t}}^{\widehat{t} + \Delta t} \hspace{-0.5cm} v(\cdot,s) \, ds \ - \ v(\cdot,\widehat{t}) \biggr) {\biggr\|}_{L^{q}(E)} 
\end{aligned}
\end{equation*}
\vspace{-0.7cm}
\begin{equation*}
\begin{aligned}
 \displaystyle
 & \! \! \leq
 \frac{1}{h} \ \limsup_{\Delta t \to 0} \biggl( \Bigl\| \frac{1}{\Delta t} \int_{\widehat{t} + h}^{\widehat{t} + \Delta t + h} \hspace{-0.5cm} v(\cdot,s) \, ds - \int_{\widehat{t} + h}^{\widehat{t} + \Delta t + h} \hspace{-0.5cm} w(\cdot,s) \, ds + \int_{\widehat{t} + h}^{\widehat{t} + \Delta t + h} \hspace{-0.5cm} w(\cdot,s) \, ds
 \qquad \qquad \\
 & \qquad \qquad \qquad \qquad \qquad \qquad \qquad - w(\cdot,\widehat{t}+h) + w(\cdot,\widehat{t}+h) - v(\cdot,\widehat{t} + h) {\Bigr\|}_{L^{q}(E)} \biggr) \\
 & \qquad \quad + \displaystyle \ \frac{1}{h} \ \limsup_{\Delta t \to 0} \biggl(\Bigl\| \frac{1}{\Delta t} \int_{\widehat{t}}^{\widehat{t} + \Delta t} \hspace{-0.5cm} v(\cdot,s) \, ds - \int_{\widehat{t}}^{\widehat{t} + \Delta t} \hspace{-0.5cm} w(\cdot,s) \, ds + \int_{\widehat{t}}^{\widehat{t} + \Delta t} \hspace{-0.5cm} w(\cdot,s) \, ds\ \\
 & \qquad \qquad \qquad \qquad \qquad \qquad \qquad \qquad \qquad \qquad - w(\cdot,\widehat{t}) + w(\cdot,\widehat{t}) - \ v(\cdot,\widehat{t}) {\Bigr\|}_{L^{q}(E)} \biggr)
\end{aligned}
\end{equation*}
\vspace{-0.5cm}
\begin{equation*}
\begin{aligned}
 \displaystyle
 & \quad \ \leq \biggl( \frac{1}{h} \ \limsup_{\Delta t \to 0} \frac{1}{|\Delta t|} \ \Bigl| \int_{\widehat{t} + h}^{\widehat{t} + \Delta t + h} \bigl\| v(\cdot,s) - w(\cdot,s) {\bigr\|}_{L^{q}(E)} \, ds \ \Bigr|
 \qquad \qquad \qquad \qquad \qquad \qquad \qquad \\
 & \qquad \qquad \qquad \quad + \displaystyle \frac{1}{h} \ \limsup_{\Delta t \to 0} \Bigl\| \frac{1}{\Delta t} \int_{\widehat{t} + h}^{\widehat{t} + \Delta t + h} \bigl(w(\cdot,s) - w(\cdot,\widehat{t}+h) \bigr) \, ds {\Bigr\|}_{L^{q}(E)} \\
 & \qquad \qquad \qquad \qquad \qquad \qquad \qquad + \ \displaystyle \frac{1}{h} \ \limsup_{\Delta t \to 0} \bigl\| w(\cdot,\widehat{t}+h) - v(\cdot,\widehat{t}+h) {\bigr\|}_{L^{q}(E)} \biggr) \\
 & \qquad \quad \displaystyle + \ \biggl( \frac{1}{h} \ \limsup_{\Delta t \to 0} \frac{1}{|\Delta t|} \ \Bigl| \int_{\widehat{t}}^{\widehat{t} + \Delta t} \bigl\| v(\cdot,s) - w(\cdot,s) {\bigr\|}_{L^{q}(E)} \, ds \Bigr| \\
 & \qquad \qquad \qquad \qquad \qquad \quad + \ \displaystyle \frac{1}{h} \ \limsup_{\Delta t \to 0} \Bigl\| \frac{1}{\Delta t} \int_{\widehat{t}}^{\widehat{t} + \Delta t} \bigl(w(\cdot,s) - w(\cdot,\widehat{t}) \bigr) \, ds {\Bigr\|}_{L^{q}(E)}  \\
 & \qquad \qquad \qquad \qquad \qquad \qquad \qquad \qquad \qquad \ + \ \displaystyle \frac{1}{h} \ \limsup_{\Delta t \to 0} \bigl\| w(\cdot,\widehat{t}) - v(\cdot,\widehat{t}) {\bigr\|}_{L^{q}(E)} \biggr) \qquad
\end{aligned}
\end{equation*}
\vspace{-0.5cm}
\begin{equation*}
\begin{aligned}
 \displaystyle
 & \ \ \ \leq \biggl( \frac{1}{h} \ \limsup_{\Delta t \to 0} \frac{1}{\Delta t} \ \int_{\widehat{t} + h}^{\widehat{t} + \Delta t + h} \bigl\| v(\cdot,s) - w(\cdot,s) {\bigr\|}_{L^{q}(E)} \, ds \ + \ 0 \ + \ \varepsilon \ \biggr) \qquad \qquad \qquad \qquad \qquad \\
 & \qquad \qquad \displaystyle + \ \biggl( \frac{1}{h} \ \limsup_{\Delta t \to 0} \frac{1}{\Delta t} \ \int_{\widehat{t}}^{\widehat{t} + \Delta t} \bigl\| v(\cdot,s) - w(\cdot,s) {\bigr\|}_{L^{q}(E)} \, ds \ + \ 0 \ + \ \varepsilon \biggr).
\end{aligned}
\end{equation*}
So that we have
\begin{equation}\tag{4.3$g$}\label{4.3g}
\begin{aligned}
 \displaystyle
 & \limsup_{\Delta t \to 0} {\Bigl\| \frac{1}{\Delta t} \Bigl( {v}_{h}(\cdot,\widehat{t} + \Delta t) - {v}_{h}(\cdot, \widehat{t}) \Bigr) - \frac{1}{h} \Bigl( v(\cdot,\widehat{t}+h) - v(\cdot,\widehat{t}) \Bigr) \Bigr\|}_{L^{q}(E)} \leq \\
 & \qquad \qquad
 \leq \displaystyle \frac{1}{h} \ \limsup_{\Delta t \to 0} \frac{1}{\Delta t} \ \int_{\widehat{t} + h}^{\widehat{t} + \Delta t + h} \bigl\| v(\cdot,s) - w(\cdot,s) {\bigr\|}_{L^{q}(E)} \, ds \  + \\
 & \qquad \qquad \qquad \qquad
 +\ \frac{1}{h} \ \limsup_{\Delta t \to 0} \frac{1}{\Delta t} \ \int_{\widehat{t}}^{\widehat{t} + \Delta t} \bigl\| v(\cdot,s) - w(\cdot,s) {\bigr\|}_{L^{q}(E)} \, ds \  + \frac{2}{h} \ \varepsilon.
\end{aligned}
\end{equation}

Now,
\begin{equation}\tag{4.3$h$}\label{4.3h}
 \begin{array}{rcl}
  \displaystyle
  \frac{1}{|\Delta t|} \biggl| \int_{\widehat{t}+h}^{\widehat{t}+\Delta t+h} \hspace{-0.8cm} \| v(\cdot,s) - w(\cdot,s) {\|}_{L^{q}(E)} ds \ \biggr| \! & \leq &
  \! \displaystyle \frac{1}{|\Delta t|} \biggl| \int_{\widehat{t}+h}^{\widehat{t}+\Delta t+h} \hspace{-0.8cm} \| v(\cdot,s) - w(\cdot,s) {\|}_{L^{q}(E)}^{r} ds {\biggr|}^{\frac{1}{r}} \ {|\Delta t|}^{1 - \frac{1}{r}} \\
  & = & \! \displaystyle {\Biggl( \frac{1}{\Delta t} \int_{\widehat{t} + h}^{\widehat{t} + \Delta t + h} \bigl\| v(\cdot,s) - w(\cdot,s) {\bigr\|}_{L^{q}(E)}^{r} ds \Biggr)}^{\frac{1}{r}}
 \end{array}
\end{equation}
and, similarly,
\begin{equation}\tag{4.3$h$\,'}\label{4.3h'}
 \displaystyle
 \frac{1}{|\Delta t|} \biggl| \int_{\widehat{t}}^{\widehat{t}+\Delta t} \bigl\| v(\cdot,s) - w(\cdot,s) {\bigr\|}_{L^{q}(E)} ds \ \biggr| \leq {\Biggl( \frac{1}{\Delta t} \int_{\widehat{t}}^{\widehat{t} + \Delta t} \bigl\| v(\cdot,s) - w(\cdot,s) {\bigr\|}_{L^{q}(E)}^{r} ds \Biggr)}^{\frac{1}{r}}.
\end{equation}

By (\ref{4.3d}) and (\ref{4.3e}), $\bigl($since $\widehat{t} \in {\ring{I}}_{h} \backslash {Z}_{**} \bigr)$ we have
\begin{equation*}
 \displaystyle
 \lim_{\Delta t \to 0} \frac{1}{\Delta t} \int_{\widehat{t}+h}^{\widehat{t}+\Delta t+h} \bigl\| v(\cdot,s) - w(\cdot,s) {\bigr\|}_{L^{q}(E)}^{r} ds
 = \bigl\| v(\cdot,\widehat{t}+h) - w(\cdot,\widehat{t}+h) {\bigr\|}_{L^{q}(E)}^{r}
 \leq {\varepsilon}^{r}
\end{equation*}
and
\begin{equation*}
 \displaystyle
 \lim_{\Delta t \to 0} \frac{1}{\Delta t} \int_{\widehat{t}}^{\widehat{t}+\Delta t} \bigl\| v(\cdot,s) - w(\cdot,s) {\bigr\|}_{L^{q}(E)}^{r} ds
 = \bigl\| v(\cdot,\widehat{t}) - w(\cdot,\widehat{t}) {\bigr\|}_{L^{q}(E)}^{r}
 \leq {\varepsilon}^{r},
\end{equation*}
so that, by (\ref{4.3h}) e (\ref{4.3h'}), we have
\begin{equation*}
 \displaystyle
 \limsup_{\Delta t \to 0} \biggl| \frac{1}{\Delta t} \int_{\widehat{t}+h}^{\widehat{t}+\Delta t+h} \bigl\| v(\cdot,s) - w(\cdot,s) {\bigr\|}_{L^{q}(E)} ds \ \biggr| \ \leq \ \varepsilon
\end{equation*}
and
\begin{equation*}
 \displaystyle
 \limsup_{\Delta t \to 0} \biggl| \frac{1}{\Delta t} \int_{\widehat{t}}^{\widehat{t}+\Delta t} \bigl\| v(\cdot,s) - w(\cdot,s) {\bigr\|}_{L^{q}(E)} ds \ \biggr| \ \leq \ \varepsilon.
\end{equation*}

\vspace{0.1cm}
Therefore, by (\ref{4.3g}), we obtain
\begin{equation*}
 \begin{aligned}
  \displaystyle
  & \limsup_{\Delta t \to 0} \Bigl\| \frac{1}{\Delta t} \Bigl( {v}_{h}(\cdot,\widehat{t} + \Delta t) - {v}_{h}(\cdot, \widehat{t}) \Bigr) - \frac{\bigl( v(\cdot,\widehat{t}+h) - v(\cdot,\widehat{t}) \bigr)}{h} {\Bigr\|}_{L^{q}(E)} \leq \\
  & \qquad \qquad \leq \frac{1}{h} \limsup_{\Delta t \to 0} \frac{1}{\Delta t} \ \int_{\widehat{t} + h}^{\widehat{t} + \Delta t + h} \bigl\| v(\cdot,s) - w(\cdot,s) {\bigr\|}_{L^{q}(E)} \, ds \\
  & \qquad \qquad \qquad \qquad + \ \frac{1}{h} \ \limsup_{\Delta t \to 0} \frac{1}{\Delta t} \ \int_{\widehat{t}}^{\widehat{t} + \Delta t} \bigl\| v(\cdot,s) - w(\cdot,s) {\bigr\|}_{L^{q}(E)} \, ds \ \\
  & \qquad \qquad \qquad \qquad \qquad \qquad + \frac{2}{h} \varepsilon \\
  & \qquad \qquad \leq \frac{1}{h} \ \varepsilon + \frac{1}{h} \ \varepsilon + \frac{2}{h} \ \varepsilon \ = \ \frac{4}{h} \ \varepsilon.
 \end{aligned}
\end{equation*}
Because $\varepsilon > 0$ is arbitrary (and $h>0$ is fixed), this shows (\ref{4.3f}), and the proof of Lemma \ref{4.3} is now complete.

\vspace{-0.8cm}
\begin{flushright}
(Lemma \ref{4.3}) \ $\square$
\end{flushright}

%
%

\section{About the integration of the Steklov averages}\label{5}

\vspace{0.3cm}
From Lemma \ref{4.3}, if
$v(\cdot,t) \ \in \ L^{r}_{loc}\bigl(I,L^{q}(E)\bigr)$
(where $1 \leq r < \infty$, $1 \leq q \leq \infty$), we always obtain that
\begin{equation*}
 \displaystyle
 {v}_{h}(\cdot,t) = \frac{1}{h} \int_{t}^{t+h} v(\cdot,s) \, ds \qquad \Bigl(\in C^{0}\bigl({\ring{I}}_{h},L^{q}(E)\bigr)\Bigr)
\end{equation*}
is (strongly) differentiable $\bigl($as a map from ${\ring{I}}_{h}$ to $L^{q}(E) \bigr)$ at almost every point $t \in {\ring{I}}_{h}$, with
\begin{equation*}
 \displaystyle {\bigl({v}_{h}\bigr)}_{t}(\cdot,t) = \frac{v(\cdot,t+h) - v(\cdot,t)}{h}, \ \ \forall \ t \in {\ring{I}}_{h} \backslash {Z},
\end{equation*}
for some ${Z} \subseteq I$, with $|{Z}|=0$, since
\begin{equation*}
 \displaystyle
 \frac{v(\cdot,t+h) - v(\cdot,t)}{h} \in L^{1}_{loc}\bigl({\ring{I}}_{h},L^{q}(E)\bigr).
\end{equation*}

This gives that, for any compact interval $[{t}_{1},{t}_{2}] \subseteq {\ring{I}}_{h}$, we have:
\begin{equation*}
 \begin{aligned}
  \displaystyle
  \int_{t_{1}}^{t_{2}} {\bigl({v}_{h}\bigr)}_{t}(\cdot,t) \, dt
  & = \int_{t_{1}}^{t_{2}} \frac{v(\cdot,t+h) - v(\cdot,t)}{h} \, dt \\
  & = \frac{1}{h} \int_{t_{1}}^{t_{2}} v(\cdot,t+h) \, dt - \frac{1}{h} \int_{t_{1}}^{t_{2}} v(\cdot,t) \, dt \\
  & = \frac{1}{h} \int_{t_{1}+h}^{t_{2}+h} v(\cdot,s) \, ds - \frac{1}{h} \int_{t_{1}}^{t_{2}} v(\cdot,s) \, ds \\
  & = \frac{1}{h} \int_{t_{2}}^{t_{2}+h} v(\cdot,s) \, ds - \frac{1}{h} \int_{t_{1}}^{t_{1}+h} v(\cdot,s) \, ds \\
  & = {v}_{h}(\cdot,{t}_{2}) - {v}_{h}(\cdot,{t}_{1}), \ \ \forall \ [{t}_{1},{t}_{2}] \subseteq {\ring{I}}_{h},
 \end{aligned}
\end{equation*}
i.e.,
\begin{equation*}
 \displaystyle
 \int_{t_{1}}^{t_{2}} {\bigl({v}_{h}\bigr)}_{t}(\cdot,t) \ dt
 = {v}_{h}(\cdot,{t}_{2}) - {v}_{h}(\cdot,{t}_{1}), \ \ \forall \ {t}_{1} \leq {t}_{2} \in {\ring{I}}_{h}.
\end{equation*}

\vspace{0.3cm}
This is a special case of the followings more general results: Lemma \ref{5.1} and Lemma \ref{5.2} above.

%
%

\vspace{0.3cm}
\begin{lem}\label{5.1}
(Fundamental Theorem of Calculus - version 1)

Given an interval $I\subseteq \R$, $E \in \mathcal{M}(\R^{n})$, $1\leq q\leq\infty$, $t_{0} \in \ring{I}$ and $F_{0} \in L^{q}(E)$, let
\begin{equation*}
 \displaystyle
 f(\cdot, t) \in C^{0}(\ring{I}, L^{q}(E))
\end{equation*}
and
\begin{equation*}
 \displaystyle
 F(\cdot,t) = F_{0}+ \int_{t_{0}}^{t}f(\cdot, s) \, ds, \ \ \forall \ t \in \ring{I}.
\end{equation*}
Then
\begin{equation}\tag{5.1$a$}\label{5.1a}
 \displaystyle
 F_{t}(\cdot,t) = f(\cdot,t), \ \ \forall \ t \in \ring{I}
\end{equation}
and
\begin{equation}\tag{5.1$a$\,'}\label{5.1a'}
 \displaystyle
 \int_{t_{1}}^{t_{2}} f(\cdot,t) \, dt = F(\cdot,t_{2}) - F(\cdot,t_{1}), \quad \forall \ t_{1}\leq t_{2} \in \ring{I}.
\end{equation}
\end{lem}

\vspace{0.3cm} \noindent
{\bf Proof.}
(\ref{5.1a}) follows by direct computation of
\begin{equation*}
 \displaystyle
 F_{t}(\cdot,t) = \lim_{\triangle t\to0}\frac{1}{\triangle t} \left[ F(\cdot,t+\triangle t) - F(\cdot,t) \right]
\end{equation*}
and of the continuity of $f(\cdot, t)$, since, for each $t\in \ring{I}$, we have:
\begin{equation*}
 \displaystyle
 F_{t}(\cdot,t) = \lim_{\triangle t \to 0}
 \frac{1}{\triangle t} \left[ F(\cdot,t+\triangle t) - F(\cdot,t) \right] = \lim_{\triangle t \to 0} \frac{1}{\triangle t} \int_{t}^{t+\triangle t} f(\cdot,s) \, ds
 = f(\cdot,t).
\end{equation*}

On the other hand, (\ref{5.1a'}) is trivial:
\begin{equation*}
 \displaystyle \int_{t_1}^{t_2} f(\cdot,t) \, dt = \int_{t_0}^{t_2} f(\cdot,t) \, dt - \int_{t_0}^{t_1} f(\cdot,t) dt = F(\cdot,t_2) - F(\cdot, t_1).
\end{equation*}

\vspace{-0.8cm}
\begin{flushright}
(Lemma \ref{5.1}) \ $\square$
\end{flushright}

%
%

\vspace{0.3cm}
The next extension of Lemma \ref{5.1}, when $f(\cdot,t)\in L^{1}_{loc}(\ring{I},L^q(E))$
(for example, if $f(\cdot,t)\in L^{r}_{loc}(\ring{I},L^q(E)$, for some $1\leq r\leq \infty$), is also worth mentioning.

\vspace{0.3cm}
\begin{lem}\label{5.2}
(Fundamental Theorem of Calculus - version 2)

Let $I\subset \mathbb{R}$ an interval, $E\in \mathcal{M}(\mathbb{R}^n)$, $1\leq q\leq \infty$, $t_0\in\ring{I}$ and $F_0\in L^q(E)$. Given
\begin{equation}\tag{5.2$a$}\label{5.2a}
 \displaystyle
 f(\cdot,t) \in L^{1}_{loc}(\ring I,L^{q}(E))
\end{equation}
and
\begin{equation}\tag{5.2$a$\,'}\label{5.2a'}
 \displaystyle
 F(\cdot,t) := F_{0} + \int_{t_{0}}^{t} f(\cdot,s) \, ds \quad (\forall \ t \in \ring I),
\end{equation}
then we have
\begin{equation}\tag{5.2$b$}\label{5.2b}
 \displaystyle
 F_{t}(\cdot,t) = f(\cdot,t) \ \ \text{ a.e. } t \in \ring I
\end{equation}
and
\begin{equation}\tag{5.2$b$\,'}\label{5.2b'}
 \displaystyle
 \int_{t_{1}}^{t_{2}} f(\cdot,t) \, dt = F(\cdot,t_{2}) - F(\cdot,t_{1}) \quad \forall \ t_{1}\leq t_{2} \in \ring I.
\end{equation}
\end{lem}

\vspace{0.3cm} \noindent
{\bf Remark.}
If $f(\cdot,t) \in L^{1}_{loc}(\ring I, L^{q}(E))$ and $G(\cdot,t) \in C^{0}(\ring I,L^{q}(E))$ is such that
\begin{equation}\tag{5.2$c$}\label{5.2c}
 \displaystyle
 G_{t}(\cdot,t) = f(\cdot,t), \quad \text{ a.e. } t \in \ring I,
\end{equation}
it does {\it not} follow (in general) that
\begin{equation}\tag{5.2$c$\,'}\label{5.2c'}
 \displaystyle
 \int_{t_{1}}^{t_{2}} f(\cdot,t) \, dt = G(\cdot,t_{2}) - G(\cdot,t_{1}), \ \ \text{ if } [t_{1},t_{2}] \subseteq \ring I.
\end{equation}
In fact, recall the Cantor-Lebesgue function, which already shows that the result (\ref{5.2c'}) is not valid even for real-valued $f \in L^{1}_{loc}(\ring I, \R)$.

The validity of (\ref{5.2c'}) requires that $G(\cdot,t)$ be also absolutely continuous in $\ring I$, i.e., that we have
\begin{equation*}
 \displaystyle
 G(\cdot,t) := G_{0} + \int_{t_{0}}^{t} g(\cdot,s) \, ds, \ \ \ (\forall \ t \in \ring I),
\end{equation*}
for some $G_{0} \in L^{q}(E)$ and
$g(\cdot, t) \in L^{1}(\ring{I}, L^{q}(E))$, and in this case, by (\ref{5.2b}), we have that $g(\cdot, t)= f(\cdot,t)$ a.e. $t \in \ring I$.

\vspace{0.3cm} \noindent
{\bf Proof of Lemma \ref{5.2}.}
Again, the proof of (\ref{5.2b'}) is a trivial consequence of (\ref{5.2a'}): we have
\begin{equation*}
 \displaystyle
 \int_{t_{1}}^{t_{2}} f(\cdot,t) \, dt = \int_{t_{0}}^{t_{2}} f(\cdot,t) \, dt - \int_{t_{0}}^{t_{1}} f(\cdot,t) \, dt = F(\cdot,t_{2}) - F(\cdot,t_{1}).
\end{equation*}
If $f(\cdot,t) \in C^{0}(\ring I, L^{q}(E))$, then we have $F_{t}(\cdot,t) = f(\cdot,t), \ \forall \ t \in \ring I$, (see Lemma \ref{5.1}).

In the general case where it is only assumed that $f(\cdot,t) \in L_{loc}^{1}(\ring I, L^{q}(E))$, we may proceed as follows (as in the proof of Lemma \ref{4.3}): first, let $Z_{0}\subseteq \ring I$ be a set with zero measure such that
\begin{equation*}
 \displaystyle
 f(\cdot,t) \in L^{q}(E), \ \ \forall \ t \in \ring I \! \setminus \! Z_{0}.
\end{equation*}

Taking a sequence of smooth approximations $g_{m}(\cdot,t) \in C_{c}^{0}(\ring I, L^{q}(E))$, $m=1,2,3,...,$ such that
\begin{equation}\tag{5.2$d$}\label{5.2d}
 \displaystyle
 \int_{a}^{b} {\| g_{m}(\cdot,t) - f(\cdot,t) \|}_{L^{q}(E)} \, dt \longrightarrow 0
\end{equation}
as $m \rightarrow \infty,$ for each compact $[a,b] \subseteq \ring I$, and
\begin{equation}\tag{5.2$e$}\label{5.2e}
 \displaystyle
 {\| g_{m}(\cdot,t) - f(\cdot,t) \|}_{L^{q}(E)} \longrightarrow 0
\end{equation}
as $m \rightarrow \infty$, for each $t \in \ring I \! \setminus \! Z_{*}$ (where $Z_{*}\subseteq \ring I$ is some null set with $Z_{*}\supseteq Z_{0}$);
and for each $m=1,2,,3,...$, let $Z_{m} \subseteq \ring I$ be a null set such that
\begin{equation*}
 \displaystyle
 \frac{1}{\Delta t} \int_{t}^{t+\Delta t} {\| f(\cdot,s) - g_{m}(\cdot,s) \|}_{L^{q}(E)} \, ds \longrightarrow {\| f(\cdot,t) - g_{m}(\cdot,t) \|}_{L^{q}(E)}
\end{equation*}
as $\Delta t \rightarrow 0$, for each $t \in \ring I \! \setminus \! Z_{m}$ (by the standard Lebesgue's differentiation theorem on $L^{1}_{loc}(\ring I, \R)$).
Then, taking
\begin{equation}\tag{5.2$f$}\label{5.2f}
 \displaystyle
 Z_{**} := Z_{0} \bigcup Z_{*} \bigcup \left( \bigcup_{m=1}^{\infty} Z_{m} \right),
\end{equation}
we have $Z_{**}\subseteq \ring I$, $Z_{**} \supseteq Z_{0}$, $|Z_{**}|=0$, and
\begin{equation}\tag{5.2$g$}\label{5.2g}
 \displaystyle
 \lim_{\Delta t \rightarrow0} {\left\| \frac{1}{\Delta t} \bigl[F(\cdot,t+\Delta t) - F(\cdot,t)\bigr] - f(\cdot,t) \right\|}_{L^{q}(E)} = 0
\end{equation}
for every $t \in \ring I \! \setminus \! Z_{**}$, thus showing (\ref{5.2b}).

To finish the proof of Lemma \ref{5.2}, it remains to shows (\ref{5.2g}).

\noindent
{\it Claim:} (\ref{5.2g}) is true.

Indeed, given $\hat{t} \in \ring I \! \setminus \! Z_{**}$ $\bigl($where $Z_{**} \subseteq \ring I$ is given in (\ref{5.2f})$\bigr)$, we now show that (\ref{5.2g}) holds at $t=\hat{t}$, i.e., we have
\begin{equation}\tag{5.2$h$}\label{5.2h}
 \displaystyle
 \limsup_{\Delta t \rightarrow 0} {\left\| \frac{1}{\Delta t} \big[ F(\cdot,\hat{t}+\Delta t) - F(\cdot,\hat{t}\,) \big] - f(\cdot,\hat{t}\,) \right\|}_{L^{q}(E)} = 0.
\end{equation}
In fact, given $\varepsilon>0$, let $g \equiv g_{m}$ be an approximant in the sequence $(g_{m})_{m}$ such that, as $m \to \infty$, we have
\begin{equation}\tag{5.2$h$\,'}\label{5.2h'}
 \displaystyle
 {\|f(\cdot,\hat{t}\,) - g(\cdot,\hat{t}\,) \|}_{L^{q}(E)} \leq \varepsilon.
\end{equation}
(\ref{5.2h'}) comes from the condition (\ref{5.2e}) (the property (\ref{5.2d}) will not be used here).

Writting (for $|\Delta t|$ small, namely, $0 < |\Delta t| \leq \hat{\delta}$, where $\hat{\delta} > 0$ is such that $[\hat{t} - \hat{\delta},\hat{t} - \hat{\delta}] \subseteq \ring I$):
\begin{equation*}
 \begin{aligned}
  \displaystyle
  & \frac{1}{\Delta t} \big[F(\cdot,\hat{t}+\Delta t)-F(\cdot, \hat{t}) \big]-f(\cdot, \hat{t}) = \frac{1}{\Delta t} \int_{\hat{t}}^{\hat{t}+ \Delta t}f(\cdot,s)ds -f(\cdot,\hat{t}) \\
  & \qquad = \frac{1}{\Delta t} \int_{\hat{t}}^{\hat{t} + \Delta t} \big[ f(\cdot,s) - g(\cdot,s) \big] \, ds + \frac{1}{\Delta t} \int_{\hat{t}}^{\hat{t} + \Delta t} g(\cdot,s) \, ds - g(\cdot,\hat{t}\,) + g(\cdot,\hat{t}\,)- f(\cdot,\hat{t}\,),
 \end{aligned}
\end{equation*}
we get
\begin{equation*}
 \begin{aligned}
  \displaystyle
  & {\left\| \frac{1}{\Delta t} \big[ F(\cdot,\hat{t}+\Delta t) - F(\cdot,\hat{t})  \big] - f(\cdot,\hat{t}\,) \right\|}_{L^{q}(E)} \leq \\
  & \qquad \leq \frac{1}{\Delta t} \int_{\hat{t}}^{\hat{t}+\Delta t} {\| f(\cdot,s) - g(\cdot,s) \|}_{L^{q}(E)} \, ds + {\left\| \frac{1}{\Delta t} \int_{\hat{t}}^{\hat{t}+\Delta t} g(\cdot,s) \, ds - g(\cdot,\hat{t}\,) \right\|}_{L^{q}(E)} \\
  & \qquad \qquad \qquad \qquad \qquad \qquad \qquad \qquad \qquad \qquad \qquad \qquad \qquad
  + {\left\| g(\cdot,\hat{t}\,) - f(\cdot,\hat{t}\,) \right\|}_{L^{q}(E)} \\
  & \qquad \leq \frac{1}{\Delta t} \int_{\hat{t}}^{\hat{t}+\Delta t} {\| f(\cdot,s) - g(\cdot,s) \|}_{L^{q}(E)} \, ds + \frac{1}{\Delta t} \int_{\hat{t}}^{\hat{t}+\Delta t} {\left\| g(\cdot,s) - g(\cdot,\hat{t}\,) \right\|}_{L^{q}(E)} \, ds + \varepsilon,
 \end{aligned}
\end{equation*}
for all $\Delta t \in \R$ with $0<|\Delta t|\leq\hat{t}$, where $\hat{t}>0$ is such that $[\hat{t}-\hat{\delta},\hat{t}+\hat{\delta}] \subseteq \ring I$. This gives
\begin{equation*}
 \begin{aligned}
  \displaystyle
  & \limsup_{\Delta t \rightarrow 0} {\left\| \frac{1}{\Delta t} \big[ F(\cdot,\hat{t}+\Delta t) - F(\cdot, \hat{t}) \big] - f(\cdot,\hat{t}\,) \right\|}_{L^{q}(E)} \leq \\
  & \quad \leq \limsup_{\Delta t \rightarrow 0} \left( \! \frac{1}{\Delta t} \! \int_{\hat{t}}^{\hat{t}+\Delta t} \!\!\! {\| f(\cdot,s) - g(\cdot,s) \|}_{L^{q}(E)} \! ds + \frac{1}{\Delta t} \! \int_{\hat{t}}^{\hat{t}+\Delta t} \!\!\! {\| g(\cdot,s) - g(\cdot,\hat{t}\,) \|}_{L^{q}(E)} \! ds \! \right) \! + \varepsilon \\
  & \quad = \limsup_{\Delta t \rightarrow 0} \frac{1}{\Delta t} \int_{\hat{t}}^{\hat{t}+\Delta t} {\| f(\cdot,s) - g(\cdot,s) \|}_{L^{q}(E)} \, ds + \varepsilon \\
  & \quad = {\| f(\cdot,\hat{t}\,) - g(\cdot,\hat{t}\,) \|}_{L^{q}(E)}
  + \varepsilon
  \ \leq \ \varepsilon +\varepsilon
 \end{aligned}
\end{equation*}
i.e., we have
\begin{equation}\tag{5.2$i$}\label{5.2i}
 \displaystyle
 \limsup_{\Delta t \rightarrow 0} {\left\| \frac{1}{\Delta t} \big[ F(\cdot,\hat{t}+\Delta t) - F(\cdot,\hat{t}\,) \big] - f(\cdot,\hat{t}\,) \right\|}_{L^{q}(E)} \leq 2 \varepsilon,
\end{equation}
for any $\varepsilon>0.$

Because $\varepsilon > 0$ in (\ref{5.2i}) is arbitrary, this shows (\ref{5.2h}), where $t = \hat{t} \in \ring I \! \setminus \! Z_{**}$ is also arbitrary. This completes the proof of Lemma \ref{5.2}.

\vspace{-0.5cm}
\begin{flushright}
(Lemma \ref{5.2}) \ $\square$
\end{flushright}

\noindent
{\bf Remark.}
Lemma \ref{5.2} can be used to give a shorter proof to Lemma \ref{4.3}, since (picking $t_{0} \in \ring I$ arbitrary) we have
\begin{equation*}
 \displaystyle
 v_{h}(\cdot,t) = \frac{1}{h} \int_{t}^{t+h} v(\cdot,s) \, ds = \frac{1}{h} \int_{t_{0}}^{t+h} v(\cdot,s) \, ds - \frac{1}{h} \int_{t_{0}}^{t} v(\cdot,s) \, ds,
\end{equation*}
and, by Lemma \ref{5.2}, since $v(\cdot,t) \in L^{r}(I,L^{q}(E)) \subseteq L_{loc}^{1}(\ring I, L^{q}(E))$, we have
\begin{equation*}
 \displaystyle
 \frac{\partial}{\partial t} \int_{t_{0}}^{t+h} v(\cdot,s) \, ds = v(\cdot,t+h),
 \ \ \text{ a.e. } t \in \ring I_{h},
\end{equation*}
(by a similar argument to that given in the proof of Lemma \ref{5.2}); and
\begin{equation*}
 \displaystyle
 \frac{\partial}{\partial t} \int_{t_{0}}^{t} v(\cdot,s) \, ds = v(\cdot,t),
 \ \ \text{ a.e. } t \in \ring I.
\end{equation*}

%
%

\vspace{0.3cm}
\begin{lem}\label{5.3}
(Integration by parts; special case: $F,G(\cdot,t) \in C^{1}(\ring I,L^{q}(E))$)

Let $I\subseteq \R$ (an interval), $E \in \mathcal M(\R^{n})$,
$1 \leq q \leq \infty$, $t_0,t_1 \in \ring I$ and $F_0, G_1 \in L^q(E)$. Given
\begin{equation}\tag{5.3$a$}\label{5.3a}
 \displaystyle
 f(\cdot,t), \ g(\cdot,t) \in C^0(\ring I,L^q(E)),
\end{equation}
\begin{equation*}
 \displaystyle
 F(\cdot,t) := F_0 + \int_{t_0}^{t} f(\cdot,s) \, ds, \ \ \forall \ t \in \ring{I},
\end{equation*}
and
\begin{equation*}
 \displaystyle
 G(\cdot,t) := G_1 + \int_{t_1}^{t} g(\cdot,s) \, ds, \ \ \forall \ t \in \ring{I},
\end{equation*}
then, for any compact interval $[a,b] \subseteq \ring I$, we have:
\begin{equation}\tag{5.3$b$}\label{5.3b}
 \displaystyle
 \int_{a}^{b} f(\cdot,t) G(\cdot,t) \, dt
 = F(\cdot,b) G(\cdot,b) - F(\cdot,a) G(\cdot,a)
 - \int_{a}^{b} F(\cdot,t) g(\cdot,t) \, dt,
\end{equation}
i.e.,
\begin{equation}\tag{5.3$b$\,'}\label{5.3b'}
 \displaystyle
 \int_{a}^{b} F_{t}(\cdot,t) G(\cdot,t) \, dt
 = F(\cdot,b) G(\cdot,b) - F(\cdot,a) G(\cdot,a)
 - \int_{a}^{b} F(\cdot,t) G_t(\cdot,t) \, dt.
\end{equation}
\end{lem}

%
%
%
%

\vspace{0.3cm} \noindent
{\bf Proof.}
Given $[a,b]\subset \ring I$, (using the Lebesgue's dominated convergence theorem) we have
\begin{equation*}
 \begin{aligned}
  \displaystyle
  & \! \! \! \! \! \int_{a}^{b} f(\cdot,t)G(\cdot,t) \, dt = \int_{a}^{b} \left[ \lim_{\Delta t \to 0} \frac{1}{\Delta t} \int_t^{t+\Delta t} f(\cdot,s) \, ds \right] \cdot G(\cdot,t) \, dt \\
  & \ \ = \lim_{\Delta t\to 0} \int_{a}^{b} \frac{1}{\Delta t} \left[ \int_t^{t+\Delta t} f(\cdot,s) \, ds \right] \cdot G(\cdot,t) \, dt \\
  & \ \ = \lim_{\Delta t\to 0} \frac{1}{\Delta t} \int_{a}^{b} \left[ F(\cdot,t+\Delta t) - F(\cdot,t) \right] G(\cdot,t) dt \\
  & \ \ = \lim_{\Delta t\to 0} \frac{1}{\Delta t} \left[ \int_{a}^{b} F(\cdot,t+\Delta t) G(\cdot,t) \, dt - \int_{a}^{b} F(\cdot,t) G(\cdot,t) \, dt \right] \\
  & \ \ = \lim_{\Delta t\to 0} \frac{1}{\Delta t} \left[ \int_{a}^{b} F(\cdot,t) G(\cdot,t-\Delta t) \, dt + \int_{b}^{b+\Delta t} F(\cdot,t) G(\cdot,t-\Delta t) \, dt \right. \\
  & \qquad \qquad \qquad \qquad \displaystyle -\left. \int_{a}^{a+\Delta t} F(\cdot,t) G(\cdot,t-\Delta t) \, dt - \int_{a}^{b} F(\cdot,t) G(\cdot,t) \, dt \right]
  \qquad \quad
 \end{aligned}
\end{equation*}
\begin{equation*}
 \begin{aligned}
  \displaystyle
  & \qquad \ = \lim_{\Delta t\to 0} \left[ \int_{a}^{b} F(\cdot,t) \frac{G(\cdot,t) - G(\cdot,t-\Delta t)}{\Delta t} \, dt \right] \\
  & \qquad \ \ \ \ + \lim_{\Delta t \to 0} \frac{1}{\Delta t} \int_{b}^{b+\Delta t} F(\cdot,t) \ G(\cdot,t-\Delta t) \, dt - \lim_{\Delta t \to 0} \frac{1}{\Delta t} \int_{a}^{a+\Delta t} F(\cdot,t) \ G(\cdot,t-\Delta t) \, dt
 \end{aligned}
\end{equation*}
\begin{equation*}
 \begin{aligned}
  \displaystyle
  & \quad = - \lim_{\Delta t \to 0} \left[ \int_{a}^{b} F(\cdot,t) \ \frac{G(\cdot,t) - G(\cdot,t-\Delta t)}{\Delta t} \, dt \right] + F(\cdot,b) G(\cdot,b) - F(\cdot,a) G(\cdot,a)
 \end{aligned}
\end{equation*}
\begin{equation*}
 \begin{aligned}
  \displaystyle
  & = - \int_{a}^{b} F(\cdot,t) \left[ \lim_{\Delta t \to 0} \frac{1}{\Delta t} \int_{t-\Delta t}^{t} g(\cdot,s) \, ds \right] \, dt + F(\cdot,b) G(\cdot,b) - F(\cdot,a) G(\cdot,a) \\
  & = - \int_{a}^{b} F(\cdot,t) g(\cdot,t) \, dt + F(\cdot,b) G(\cdot,b) - F(\cdot,a) G(\cdot,a).
 \end{aligned}
\end{equation*}
This concludes the proof of Lemma \ref{5.3}.

\vspace{-0.5cm}
\begin{flushright}
(Lemma \ref{5.3}) \ $\square$
\end{flushright}

%
%

\vspace{0.3cm}
\begin{lem}\label{5.4}
(Integration by parts; general case: $F,G(\cdot,t)$ absolutely continuous in $\ring I$)

Let $I \subseteq \R$ (an interval), $E \in \mathcal M(\R^{n})$, $t_0,t_1\in \ring{I}$,
$1 \leq q \leq \infty$, and $F_0,G_1 \in L^q(E)$. Given
\begin{equation*}
 \displaystyle
 f(\cdot,t), g(\cdot,t) \in L_{loc}^{1}(\ring I,L^{q}(E)),
\end{equation*}
\begin{equation*}
 \displaystyle
 F(\cdot,t) := F_{0} + \int_{t_{0}}^{t} f(\cdot,s) \, ds \ \ (\forall \ t \in \ring I),
\end{equation*}
and
\begin{equation*}
 \displaystyle
 G(\cdot,t) := G_{1} + \int_{t_{1}}^{t} g(\cdot,s) \, ds \ \ (\forall \ t \in \ring I),
\end{equation*}
then, for any compact $[a,b] \subseteq \ring I$, we have
\begin{equation}\tag{5.4$a$}\label{5.4a}
 \displaystyle
 \int_{a}^{b} f(\cdot,t) G(\cdot,t) \, dt = F(\cdot,b) G(\cdot,b) - F(\cdot,a) G(\cdot,a) - \int_{a}^{b} F(\cdot,t) g(\cdot,t) \, dt,
\end{equation}
i.e.,
\begin{equation}\tag{5.4$a$\,'}\label{5.4a'}
 \displaystyle
 \int_{a}^{b} F_{t}(\cdot,t) G(\cdot,t) \, dt = F(\cdot,b) G(\cdot,b) - F(\cdot,a) G(\cdot,a) - \int_{a}^{b} F(\cdot,t) G_{t}(\cdot,t) \, dt.
\end{equation}
\end{lem}

\vspace{0.3cm} \noindent
{\bf Proof.}
Taking $\big(f_{m}(\cdot,t) \big)_{m}$ and $\big(g_{m}(\cdot,t) \big)_{m}$, with $f_{m}(\cdot,t), \, g_{m}(\cdot,t) \in C_{c}^{0}(\ring I, L^{q}(E))$, such that
\begin{equation}\tag{5.4$b$}\label{5.4b}
 \displaystyle
 \int_{\alpha}^{\beta} {\|f_{m}(\cdot,t) - f(\cdot,t)\|}_{L^{q}(E)} \, dt
 \rightarrow 0, \ \ \int_{\alpha}^{\beta} {\| g_{m}(\cdot,t) - g(\cdot,t) \|}_{L^{q}(E)} \, dt \longrightarrow 0, \ \ \text{ as } m \to \infty,
\end{equation}
for each compact interval $[\alpha,\beta] \in \ring I$ we set $F_{m}(\cdot,t)$,
$G_{m}(\cdot,t) \in C(\ring I, L^{q}(E))$ given by
\begin{equation}\tag{5.4$b$\,'}\label{5.4b'}
 \displaystyle
 F_{m}(\cdot,t):= F_{0} + \int_{t_{0}}^{t}f_{m}(\cdot,s) \, ds, \ \
 G_{m}(\cdot,t):= G_{1} + \int_{t_{1}}^{t}g_{m}(\cdot,s) \, ds \ \ (\forall \ t \in \ring I),
\end{equation}
so that we have $F_{m}(\cdot,t) \rightarrow F(\cdot,t), \ G_{m}(\cdot,t) \rightarrow G(\cdot,t)$ in $L^{q}(E)$, uniformly in $t \in [\alpha,\beta]$ (for any compact $[\alpha,\beta] \subseteq \ring I$), as $m\rightarrow\infty$.

By Lemma \ref{5.3} (given any compact $[a,b] \subseteq \ring I$), $\forall \ m$ we have
\begin{equation}\tag{5.4$d$\,''}\label{5.4d''}
 \begin{aligned}
  \displaystyle
  \int_{a}^{b} f_{m}(\cdot,t) G_{m}(\cdot,t) \, dt
     & = F_{m}(\cdot,b)G_{m}(\cdot,b) \\
     & \displaystyle - F_{m}(\cdot,a) G_{m}(\cdot,a) - \int_{a}^{b} F_{m}(\cdot,t) g_{m}(\cdot,t) \, dt.
 \end{aligned}
\end{equation}
Letting $m \rightarrow \infty$ in (\ref{5.4b'}), we obtain (\ref{5.4a}).

\vspace{-0.5cm}
\begin{flushright}
(Lemma \ref{5.4}) \ $\square$
\end{flushright}

%
%

\vspace{0.3cm} \noindent
{\bf Acknowledgements.}
This work was partially supported by CNPq (Ministry of Science and Technology, Brazil), Grant \mbox{\small \#\,154037/2011-7}
and by CAPES (Ministry of Science and Technology, Brazil), Grant
\mbox{\small \#\,1212003/2013}.
The authors also express
their gratitude to Paulo R. Zingano (UFRGS, Brazil),
for some helpful suggestions
and discussions.

%
%



\renewcommand{\refname}{\normalfont\selectfont\LARGE {\bf References}}

%
%

%
\nl
\nl
\nl
{\small
\begin{minipage}[t]{10.00cm}
\mbox{\normalsize \textsc{Jocemar de Quadros Chagas}} \\
Departamento de Matem\'atica e Estat\'\i stica \\
Universidade Estadual de Ponta Grossa \\
Ponta Grossa, PR 84030-900, Brazil \\
E-mail: {\sf jocemarchagas@uepg.br} \\
\end{minipage}
\nl
\mbox{} \vspace{-0.450cm} \\
\nl
\begin{minipage}[t]{10.00cm}
\mbox{\normalsize \textsc{Nicolau Matiel Lunardi Diehl}} \\
Instituto Federal de Educa\c c\~ao, Ci\^encia e Tecnologia \\
Canoas, RS 92412-240, Brazil \\
E-mail: {\sf nicolau.diehl@canoas.ifrs.edu.br} \\
\end{minipage}
\nl
\mbox{} \vspace{-0.450cm} \\
\nl
\begin{minipage}[t]{10.00cm}
\mbox{\normalsize \textsc{Patr\'\i cia Lisandra Guidolin}} \\
Instituto Federal de Educa\c c\~ao, Ci\^encia e Tecnologia \\
Viam\^{a}o, RS 94410-970, Brazil \\
E-mail: {\sf patricia.guidolin@viamao.ifrs.edu.br} \\
\end{minipage}
\nl
\mbox{} \vspace{-0.450cm} \\
%
%
%
%

\end{document}